\input amstex
\documentstyle{amsppt}
\input xy
\xyoption{all}
\mag=\magstep2 
\pagewidth{11cm}
\NoBlackBoxes
\pageheight{15cm}
\def\a{\alpha}
\def\ol{\overline}
\def\b{\beta}
\def\gam{\gamma}
\def\Gam{\Gamma}
\def\del{\delta}
\def\lam{\lambda}
\def\ome{\omega}
\def\Ome{\Omega}

\def\A{\Cal A }
\def\E{\Cal E}

\def\G{G_2}

\def\Diff{\text{\rm Diff}}
\def\R{\Bbb R}
\def\C{\Bbb C}
\def\H{\Bbb H}

\def\CZ{\Cal Z}
\def\O{\Cal O}
\def\M{\frak M}

\def\w{\wedge}
\def\({\left(}
\def\){\right)}
\def\G{G_2}

\def\neg{\negthinspace}
\def\h{\hat}
\def\wideh{\widehat}
\def\til{\tilde}
\def\wtil{\widetilde}

\def\pa{\partial}

\def\s-style{\scriptstyle}
\def\ss-style{\scriptscriptstyle}

\def\arrow{\longrightarrow}

\def\CP{\dsize{\Bbb C} \bold P}

\document
\title\nofrills Moduli  spaces of special Lagrangians and \\
K\"ahler-Einstein structures 
\endtitle
\rightheadtext{Special Lagrangians and 
K\"ahler-Einstein structures }
\author	
Ryushi Goto
\endauthor
\affil		
Department of Mathematics,\\
 Graduate School of Science,\\
Osaka University,\\ 
\endaffil
\address
Toyonaka, Osaka, 560, Japan
\endaddress
\email
goto\@math.sci.osaka-u.ac.jp
\endemail
\abstract 
We shall construct a moduli space of pairs of K\"ahler-Einstein structures and 
special lagrangians and obtain smoothness of the moduli space of these pairs. 
Further we show that the moduli space of these pairs is locally embedded in a certain 
relative cohomology group.
\endabstract
\endtopmatter
\head 
\S0. Introduction
\endhead
Special lagrangians are calibrated submanifolds of Calabi-Yau manifolds
with  K\"ahler-Einstein structures, which  
have been extensively studied between differential geometry and
mathematical physics [8],[10],[11], [16]. 
In particular, a deformation of special
lagrangians is an intriguing topic and 
Mclean shows smoothness of the deformation space of special lagrangians, which is  
parameterized by  an open set of the first cohomology group [14].  
The moduli space of
K\"ahler-Einstein structures are constructed by Fujiki-Schumacher [4] and Tian-Todorov show
that the deformation space of K\"ahler-Einstein structures is also smooth [17],[18]. 
Special lagrangians are depending on a
choice of K\"ahler -Einstein structures. Hence it is natural to study  moduli spaces of special
lagrangians  under deformations of K\"ahler-Einstein structures.
In this paper we consider a pair consisting of a K\"ahler-Einstein structure $\Phi$ and 
a special lagrangian $M$ with respect to $\Phi$. 
These results of smoothness raise a question of 
whether we obtain a smooth moduli space of pairs of  K\"ahler-Einstein structures and special
lagrangians.
The purpose of this paper is to show smoothness of the moduli space of 
such pairs, to obtain some further properties.  The moduli space of such pairs is locally embedded
into a certain relative cohomology group.  
In order to explain our results more precisely,  we introduce our notations and
definitions. {\it A calabi-Yau manifold} is a  complex manifold with the trivial canonical line
bundle.  Let $X$ be a compact Calabi-Yau manifold with a K\"ahker form.
There exists a unique K\"ahler-Einstein form $\ome$ on $X$ with 
vanishing Ricci curvature in each K\"ahler class. We denote by $\Ome$ 
a nonzero form of type $(n,0)$, where $\dim_\C  X = n$. 
We call  a pair $\Phi =(\Ome,\ome)$ {\it a K\"ahler-Einstein structure} on $X$.
 {\it  A special lagrangian} $M$ of $X$ with respect to a K\"ahler-Einstein structure $\Phi$ is ,
by definition, a real
$n$ dimensional submanifold of $X$ satisfying equations, 
$$i^*_{\ss-style M}\Ome^{Im}=0, \quad 
i^*_{\ss-style M}\ome=0.$$ where $i_{\ss-style M}\: M \to X$ and $\Ome^{Im}$ denotes 
the imaginary part of the complex form $\Ome$.
We assume that a special lagrangian $M$ is compact.
Then we have a moduli space $\Cal P$ of pairs consisting of K\"ahler-Einstein structures and
special lagrangians ( see definition 4-1-3 in section 4). At first we show that a connected
component of the moduli space
$\Cal P$ is regarded as a certain relative moduli space $\M_{\ss-style{KE}}(X,M)$, 
where $M$ is a fixed real $n$ dimensional submanifold.
Then we obtain , 
\proclaim{Theorem } 
The relative moduli space $\M_{\ss-style{KE}}(X,M)$ is a smooth manifold. In particular,
$\M_{\ss-style{KE}}(X,M)$  is Hausdorff. 
\endproclaim
One of difficulties for a construction of the moduli space $\Cal P$ is that special lagrangians are
real calibrated manifolds, in which  we can not apply a general theory of deformations of
complex geometry, such as Kodira-Spencer theory. 
In this paper we consider $X$ as a real $2n$ dimensional manifold and 
K\"ahler-Einstein structures are geometric structures defined by  closed differential forms 
$\Phi = (\Ome,\ome)$. ( We refer to these closed differential forms  as calibrations .) 
We construct the moduli space $\M_{\ss-style {KE}}(X)$ of these calibrations by using the implicit
function theorem and obtain a smooth moduli space of K\"ahler-Einstein structures. 
Let $\Phi^0=(\Ome^0,\ome^0)$ be a K\"ahler-Einstein structure on $X$. Then 
we obtain an elliptic complex \#$_{\ss-style X}$ with respect to $\Phi^0$:
$$\CD
\Gam_{\ss-style X}(E^0_{\ss-style X} )@>d_{\ss-style X}>>
\Gam_{\ss-style X}(E^1_{\ss-style X} )@>d_{\ss-style X}>>
\Gam_{\ss-style X}(E^2_{\ss-style X} ).
\endCD
\tag \#$_{\ss-style X}$
$$
We denote by $H^i (\#_{\ss-style X})$ the cohomology group of the complex $\#_{\ss-style
X}$. Then the infinitesimal deformation at $\Phi^0$ is given by the first cohomology group 
$H^1(\#_{\ss-style X})$.
Let $i_{\ss-style M}\: M \arrow X$ be a special lagrangian with respect to $\Phi^0$. 
Then we also obtain an elliptic complex $\#_{\ss-style M}$ over $M$: 
$$
\CD 
\Gam_{\ss-style M}(E^0_{\ss-style M} )@>d_{\ss-style M}>>
\Gam_{\ss-style M}(E^1_{\ss-style M} )@>d_{\ss-style M}>>
\Gam_{\ss-style M}(E^2_{\ss-style M} ).
\endCD
\tag \#$_{\ss-style M}$
$$
By using the pull back $i^*_{\ss-style M}$, we have a surjective map of complexes, 
$\kappa\: \#_{\ss-style X} \to \#_{\ss-style M}$. Hence we have a short exact sequence of 
complexes: 
$$\CD
0@>>>\#_{ \ss-style{X,M}}@>>>
\#_{\ss-style X}@>\kappa>>\#_{\ss-style M}@>>>0,
\endCD
$$where $\#_{ \ss-style {X,M}}$ is a complex, 
$$
\CD
\Gam_{\ss-style {X,M}}(E^0_{\ss-style {X,M}} )@>d_{\ss-style {X,M}}>>
\Gam_{\ss-style {X,M}}(E^1_{\ss-style {X,M}} )@>d_{\ss-style {X,M}}>>
\Gam_{\ss-style {X,M}}(E^2_{\ss-style {X,M}} ),
\endCD
$$( see section 4 for definition).
We denote by $H^i( \#_{\ss-style {X,M}})$ the cohomology group of the complex
$\#_{\ss-style  {X,M}}$. 
 Then we have 
\proclaim{Theorem} 
Let $\M_{\ss-style {KE}}(X,M)$ be a relative moduli space as above. 
Then local coordinates at $\Phi^0$ is given by an open set of the cohomology group 
$H^1( \#_{\ss-style {X,M}})$.
\endproclaim
We also show that the cohomology group $H^1( \#_{\ss-style {X,M}})$ is a
subgroup of a relative de Rham cohomology group, which is topologically defined. 
Then we show that the moduli space $\M_{\ss-style {KE}}(X,M)$ is locally embedded
in  the relative de Rham cohomology group, which is called local Torelli type
theorem (Theorem 4-2-8).  It must be noted that the moduli space $\M_{\ss-style
{KE}}(X,M)$ is a total space of a fibre bundle, 
$$
\M_{\ss-style {KE}}(X,M)\arrow \wideh\M_0(X,M),
$$
where the base space $\wideh{\M}_0(X,M)$ is a submanifold of the moduli space of
K\"ahler-Einstein structures $\M_{\ss-style{KE}}(X)$, which corresponds to deformations 
preserving special lagrangians  (Proposition 4-2-9) and each fibre is regarded as  the moduli space
of special lagrangians with respect to a fixed K\"ahler-Einstein structure. 
Our moduli space is defined as a certain quotient by the action of the identity component of the
group of  diffeomorphisms. It is interesting to ask what is the quotient by the action of
whole diffeomorphisms. In theorem 4-2-11 we show that  such a quotient is an orbifold 
( see theorem 1-9 ).
In section 1, we discuss a
general theory of geometric structures defined by  closed differential forms and construct a
moduli space of such closed differential forms.  If a geometric structure is metrical, elliptic and
topological, we obtain  a smooth moduli space of them ( see definition 1-1,2,3). 
 Section 2 is devoted  to prove  theorems in section 1.
We show in section 3 that the K\"ahler-Einstein structure is metrical, elliptic and topological. 
Hence we obtain a smooth moduli space of K\"ahler-Einstein structures.  If we fix a class of
K\"ahler forms, then we have a smooth moduli space  of polarized Calabi-Yau structures.
In section 4 we obtain the relative moduli space $\M_{\ss-style {KE}}(X,M)$. 
In subsection 4-1 we study the cohomology groups $H^i( \#_{\ss-style X})$,$H^i(\#_{\ss-style
M})$  and $H^i(\#_{\ss-style {X,M}})$. In subsection 4-2 we construct a slice $S_{\Phi^0}(X,M)$, 
which is local coordinates of the moduli space $\M_{\ss-style {KE}}(X,M)$, and prove our main
theorem (theorem 4-1-5).  In the case of hyperK\"ahler structure, corresponding calibrated
submanifolds are holomorphic lagrangians. Our construction also holds in this case and we obtain
a smooth moduli space of pairs of hyperK\"ahler structures and holomorphic lagrangians. 
We will discuss this in a forthcoming paper. 
\head \S1. Moduli spaces of calibrations
\endhead
Let $V$ be a real vector space of dimension $n$. We denote by $\w^p V^*$ the
vector space  of $p$ forms on $V$. Let $\rho_p$ be the linear action of
$G=$GL$(V)$ on $\w^p V^*$. Then we have the action $\rho$ of $G$ on 
the direct sum
$\oplus_i
\w^{p_i} V^*$ by 
$$
\rho \: GL(V) \arrow \text{End}( \oplus_{i=1}^l\w^{p_i} V^*),
$$
$$\rho= (\rho_{p_1},\cdots,\rho_{p_l} ).
$$
We fix an element $\Phi^0 = ( \phi_1, \phi_2, \cdots ,\phi_l ) \in \oplus_i
\w^{p_i}V^*$  and consider the $G$-orbit $\O=\Cal O_{\Phi^0}$ through
$\Phi^0$: 
$$
\Cal O_{\Phi^0} = 
\{ \, \Phi = \rho_g \Phi^0 \in \oplus_i \w^{p_i} V^* \, |\, g \in G\, \}
$$
The orbit $\O_{\Phi^0}$ can be regarded as a homogeneous space,
$$
\Cal O_{\Phi^0} = G/ H,
$$where $H$ is the isotropy group 
$$
H = \{ \, g \in G \, | \, \rho_g \Phi^0 = \Phi^0 \, \}.
$$
We denote by $\Cal A (V)$ the orbit $\Cal O_{\Phi^0} = G/H$. 
The tangent space $E^1 (V)=T_{\Phi^0} \A(V)$ is given by 
$$
E^1(V)=T_{\Phi^0}\A(V) = 
\{\, \rho_\xi \Phi^0 \in \oplus_i \w^{p_i}V^* \, |\, \xi \in \frak g\, \},
$$
where $\rho$ denotes the differential representation of $\frak g$.
The vector space $E^1 (V)$ is the quotient space $\frak g/\frak h$.
We also define a vector space $E^0(V)$ by 
the interior product, 
$$\align
E^0(V) =& \{\, i_v \Phi^0 \in \oplus_i \w^{p_i-1} V^*\, |\, v \in V \, \}.
\endalign
$$
$E^2(V)$ is define as a vector space spanned by the following set, 
$$
E^2(V) =\text{Span}\{\, \theta\w \Phi \in \oplus_i \w^{p_i +1} V^* \, |\, \theta \in V^*
, \, 
\Phi \in E^1 (V) \, \}.
$$
Then we have the complex by the exterior product of  a nonzero $u \in V^*$, 
$$
\CD
E^0 (V) @>\w u>>E^1(V)@>\w u >>E^2 (V). 
\endCD
$$
\proclaim{Definition 1-1(elliptic orbits)}
An orbit $\O_{\Phi^0}$ is an elliptic orbit if the complex 
$$
\CD
E^0 (V) @>\w u>>E^1(V)@>\w u >>E^2 (V). 
\endCD
$$
is exact for any nonzero $u \in V^*.$ 
In other words, if $\a\w u =0$ for $\a\in E^1(V)$, then there exists $\b \in
E^0(V)$ such that $\a = \b \w u$.
\endproclaim
\proclaim{Definition 1-2(metrical orbits)}
Let $\O_{\ss-style\Phi^0}$ be an orbit as before.
An orbit $\O_{\Phi^0}$ is metrical if the isotropy group $H$ is a subgroup
of  O$(V)$ with respect to a metric $g_V$ on $V$.
\endproclaim
Let $X$ be a compact real manifold of dimension $n$. 
Then we define the $G/H-$bundle $\A(X)=\A_{_{\Cal O}}(X)$ by 
$$
\A_{_\O}(X) = \underset{x \in X} \to \bigcup{\A(T_x X)}\arrow X.
$$
We denote by $\E^1=\E^1(X)$ the set of $C^\infty$ global sections of $\A(X)$, 
$$
\E^1(X) = \Gam (X, \A(X)).
$$
Let $\Phi^0$ be a closed element of $\E^1$. Then we have the vector spaces $E^i
(T_xX)$  for each $x \in X$ and $i = 0,1,2$ respectively.
 We define the vector bundle $E^i_{\ss-style X}=E^i$ over $X$  as
$$
E^i_{\ss-style X}=E^i: = \underset {x \in X}\to\bigcup{E^i (T_x X)}\arrow X.
$$ for each $i = 0,1,2$.
(Note that the  fibre of $E^1$ is $\frak g/\frak h$.)
Then we have a complex $\#_{\Phi^0}$ 
$$
\CD
\Gam (E^0) @>d_0>>\Gam (E^1)@>d_1>>\Gam(E^2),
\endCD
 \tag{\#$_{\Phi^0}$}
$$
 where $\Gam (E^i)$ is the set of $C^\infty$ global sections for each vector
bundle and 
 $d_i = d|_{E^i}$ for each $i =0,1,2$.
The complex \#$_{\Phi^0}$ is a subcomplex of 
 de Rham complex:
$$
\CD 
\Gam (E^0)@>d_0>>\Gam (E^1)@>d_1>>\Gam (E^2) \\ 
@VVV                    @VVV         @VVV \\
\Gam (\oplus_i \w^{p_i-1} ) @>d>> \Gam ( \oplus_i \w ^{p_i} ) @>d>>\Gam
(\oplus_i
\w^{p_i+1} ) 
\endCD
$$
Hence there is the map $p$ from the cohomology group of the complex
\#$_{\Phi^0}$ to de Rham cohomology group: 
$$
p\:H^1 (\#_{\Phi^0} )  \arrow \underset i\to\oplus H^{p_i}(X,\R),
$$ 
where
$$
H^1 ( \#_{\Phi^0} ) = \{ \, \a \in \Gam (E^1) \, |\, d_1 \a =0 \, \} /
\{\, d\b \,| \, \b \in \Gam ( E^0)\, \}.
$$
\proclaim{Definition 1-3(Topological calibrations and topological orbits)} 
A closed element $\Phi^0 \in \E^1(X)$ is a topological calibration if the map 
$$
p\:H^1 (\#_{\Phi^0} )  \arrow \underset i\to\oplus H^{p_i}(X,\R)
$$
is injective. 
A manifold $X$ is topological with respect to an orbit $\O$ in $\oplus_i \w^{p_i }V^*$ if 
any closed element of $\E^1(X) $ is a topological calibration. 
An orbit $\Cal O$ is topological if $p$ is injective for each closed
form $\Phi^0
\in
\E^1(X)$ over any compact 
$n$ dimensional manifold $X$.
\endproclaim
\proclaim{Lemma 1-4}
Let $\Cal O$ be a metrical orbit and $\Phi^0$ an element of 
$\E^1 = \Gam (X, \A_{_\O}(X) )$.
Then there is a canonical metric $g_{\ss-style\Phi^0}$on $X$ corresponding to 
each $\Phi^0 $. 
\endproclaim
\demo{Proof}
The orbit $\O$ is defined in terms of $\Phi^0=\Phi^0_V \in \oplus_i
\w^{p_i}V^*$  on $V$. We also have $\Phi^0 (x)\in \A(T_x X)$ on each tangent
space $T_x X$.  Let $\text{Isom} (V, T_xX ) $ be the set of isomorphisms
between 
$V$ and $T_xX$.
Then define $H_x$ by 
$$
H_x = \{  h \in \text{Isom} (V, T_xX ) \, |\,  \Phi_V = h^*\Phi^0(x) \, \}.
$$
Then we see that $H_x$ is isomorphic to the isotropy group $H$. 
$h_* g_V$ defines the metric on the tangent space $T_x X$ for $h \in
H_x$. Since $H$ is a subgroup of O$(V)$, the metric $h_*g_V$ does not depend
on a choice of $h \in H_x$.
\qed\enddemo
At first we shall show the following:
\proclaim{Proposition 1-5} 
Let $\O$ be a metrical and elliptic orbit and 
$\Phi$ a topological element of $\E^1$. Then the dimension of H$^1 (\#_\Phi)$
is invariant  under deformations of $\Phi \in \Cal E^1$.
\endproclaim
\demo{Proof}
The complex $\{ \Gam (E^i), d_i \}$ is a subcomplex of de Rham complex. Hence we have 
the commutative diagram :
$$
\minCDarrowwidth{0.5cm}
\CD 
0@.0@.0\\
@VVV@VVV@VVV\\
 \Gam(E^0) @>{d_0}>>\Gam(E^1)@>{d_1}>>\Gam(E^2 )\\
 @VVV        @VVV          @VVV \\
\oplus_i \Gam (\w^{p_i-1})@>d>>\oplus_i\Gam (\w^{p_i}) @>d>>\oplus_i\Gam (\w^{p_i +1} )  \\
@VVV@VVV@VVV\\
\oplus_i \Gam (\w^{p_i-1})/\Gam(E^0)@>>>\oplus_i\Gam (\w^{p_i})/\Gam(E^1) @>>>\oplus_i\Gam (\w^{p_i +1} ) /\Gam(E^2) \\
@VVV@VVV@VVV\\
0@.0@.0
\endCD
$$ 
Let $g_\Phi$ be the metric corresponding to  the metrical calibration 
$\Phi$. 
We denote by  $E^i_\perp$ the orthogonal compliment of 
each vector bundle $E^i$ for $i = 0,1,2$. Then we have the
identification: 
$$\align
E^0_\perp &\cong \oplus_i  \w^{p_i-1}/E^0 \\
E^1_\perp &\cong \oplus_i \w^{p_i}/E^1 
\endalign
$$
We denote by $\h{E}^2$ by the image 
$$
\h{E}^2 = d_1\Gam (E^1).
$$
Then we have 
$$
\CD 
 \Gam(E^0) @>{d_0}>>\Gam(E^1)@>{d_1}>>\h{E}^2 \\
 @VVV        @VVV          @VVV \\
\oplus_i \Gam (\w^{p_i-1})@>d>>\oplus_i\Gam (\w^{p_i}) @>d>>\oplus_i\Gam (\w^{p_i +1} )  \\
@VV V@VV V@VV V\\
\Gam(E^0_\perp) @>>>\Gam(E^1_\perp) @>>>\Gam(E^2 _\perp).
\endCD
$$
Hence we have the long exact sequence. Since $\Gam(E^1) \arrow \h{E}^2$ 
is surjective and $\Phi$ is topological, we have the exact sequence: 
$$\minCDarrowwidth{ 0.7cm}\CD
0@>>> H^1 (\#)@>>>\oplus_i H^{p_i}_{\ss-style{DR}}(X) @>>> H^1 ( \#_\perp )
@>>> 0,
\endCD
\tag1 $$
When we consider symbols of differential operators in the diagram, 
we see that the complex \#$_\perp$
$$\CD
\Gam(E^0_\perp) @>>>\Gam(E^1_\perp) @>>>\Gam(E^2 _\perp)
\endCD
$$
is an elliptic complex. 
Hence H$^1(\#_\perp )$ is Kernel of the  elliptic operator.
Since \# is elliptic, H$^1(\#)$ is also kernel of the  elliptic operator. 
So the dimension of each cohomology group is an upper semi continuous.
Hence from (1) we see that the dimension of H$^1 (\#)$ is invariant 
under deformations. 
\qed\enddemo

Let $\O$ be an orbit in  $\oplus_i \w^{p_i }V^*$ . Then 
we define the moduli space $\M_{_\O} (X)$ by 
$$
\M_{_\O} (X) = \{ \, \Phi \in \E^1 \, |\, d\Phi =0 \, \} / \text{Diff}_0(X),
$$
where Diff$_0(X)$ is the identity component of the group of diffeomorphisms
for $X$.  We denote by $\widetilde{\M}_{_\O} (X)$ the set of closed elements in
$\E^1$.  We have the natural projection $\pi \: \widetilde{\M}_{_\O} (X) \to \M_{_\O}
(X)$.
We shall prove the following theorems in section 2.
\proclaim{Theorem 1-6 } 
If an orbit $\O$ is metrical, elliptic and topological, then the corresponding moduli space $\M_{_\O} (X)$ is a smooth manifold. 
( In particular $\M_{_\O} (X)$ is Hausdorff .)
Further each coordinates of $\M_{_\O} (X)$ around $\pi ( \Phi) $ is canonically
given by an open ball of the cohomology group H$^1 (\#_{\Phi}) $ for each $\Phi
\in
\til{\M}_{_\O} (X)$.
\endproclaim
Since de Rham cohomology group is invariant under the action of Diff$_0$, we have the map 
$$
P\: \M_{_\O} (X) \arrow \underset i \to \oplus H^{p_i}_{dR} (X).
$$
Then we have 
\proclaim{Theorem 1-7}
If an orbit $\O$ is metrical, elliptic and topological, then the map $P$ is locally injective. 
\endproclaim
\proclaim{Theorem 1-8} 
Let $I(\Phi$) be the isotropy group, 
$$
I(\Phi) = \{ \, f \in \text{\rm Diff}_0(X) \, |\,  f^* \Phi = \Phi \, \}.
$$
Then there is a sufficiently small  slice $S_{\Phi^0}$ at $\Phi^0 $
such that the isotropy group $I(\Phi^0)$ is a subgroup of $I(\Phi)$ for each
$\Phi\in S_{\Phi^0}$,i.e.,
$$
I(\Phi^0) \subset I(\Phi).
$$
(Our definition of the slice will be given in section 2)
\endproclaim
\proclaim{Theorem 1-9}
Let $\wtil{\M}_{\ss-style\O}(X)$ be the set of closed elements of $\E^1$. 
We denote by $\Diff(X)$ the group of diffeomorphism of $X$. There is the
action of $\Diff(X)$ on $ \wtil{\M}_{\ss-style\O}(X)$.  Then the quotient
$\wtil{\M}_{\ss-style\O}(X)/\Diff(X)$ is an orbifold.  
\endproclaim
\head 
\S2. Proof of theorems
\endhead
Let $X$ be a real $n$ dimensional compact manifold.
We denote by $C^\infty(X,\w^p)$ the set of smooth $p$ forms on $X$.
Let $L^2_s ( X, \w^p)$ be the Sobolev space and suppose that $s > k + \frac n2$., i.e., the completion 
of $C^\infty(X, \w^P)$ with respect to the Sobolev norm $\|\, \|_s$, where $k$ is
sufficiently large.  Then we have the inclusion $ L^2_s ( X, \w^p ) \arrow C^k (
X, \w^n )$.  We define $\E^1_s$ by 
$$
\E^1_s = C^k ( X, \A_\O(X) ) \cap L^2_s ( X, \oplus_{i=1}^l\w^{p_i} ).
\tag2-1$$
Then we have 
\proclaim{Lemma 2-1}
$\E^1_s$ is a Hilbert manifold. The tangent space $T_{\Phi ^0}\E^1_s$ at $\Phi^0$ is given by 
$$
T_{\Phi^0} \E^1_s = L^2_s ( X, E^1).
$$ 
\endproclaim
\demo{Proof}
We denote by exp the exponential map of Lie group $G=$GL$(n,\R)$. Then 
we have the map $k_x$ 
$$
k_x \: E^1 ( T_xX ) \arrow \A ( T_x X),
\tag2-2$$
by 
$$
k_x ( \rho_\xi\Phi^0(x) ) = \rho_{\exp \xi}\Phi^0 (x).
\tag2-3$$
for each tangent space $T_x X$.
From 2-3  we have the map $k$
$$
k \:L^2_s ( E^1) \arrow \E^1_s,
\tag2-4$$by 
$$
k|_{E^1(T_x X)} = k_x.
$$
The map  $k$ defines local coordinates of $\E^1_s$.
\qed\enddemo
For each $\Phi$ we have the vector bundles as in section one. We denote them by 
$E_{\ss-style{\Phi}}^0, E_{\ss-style{\Phi}}^1$ and $E_{\ss-style{\Phi}}^2$.
We define $\CZ$ to be the set of closed forms in $\E^1_s$,
$$
\CZ=\{ \, \Phi \in \E^1_s \, |\, d\Phi =0 \, \}.
$$
We also denote by $Z_{\Phi}$ the Hilbert space of closed forms of 
$T_\Phi \E^1_s=L_s^2 (X, E^1_{\ss-style{\Phi}})$, 
$$
Z_\Phi = \{\, a \in L_s^2 (X, E^1_{\ss-style{\Phi}})\, |\, da =0 \, \}.
$$
We then have 
\proclaim{Lemma 2-2} 
$\CZ$ is a Hilbert manifold with $T_\Phi \CZ = Z_\Phi$.
\endproclaim
\demo{Proof} 
The tangent space of the Hilbert manifold $\E^1_s$ at $\Phi$  is the Hilbert space 
$L^2_s ( X, E^1_{\ss-style{\Phi}})$. 
We then have a distribution $Z$ defined by 
the closed subspace 
$Z_\Phi$ of the tangent space $T_\Phi \E^1_s$. 
We shall show that the distribution $Z$ is integrable. 
Let $\Cal V^1$(resp.$\,\Cal V^2$) denote the Sobolev space $L^2_s(X, \oplus_i\w^{p_i})$ 
(resp. $L^2_{s-1}(X, \oplus_i\w^{p_i+1})$).
We define  a map  $\eta\: \Cal V^1 \to \Cal V^2$ by the exterior derivative
$$
a \mapsto da,
$$ where $a \in \Cal V^1$. 
Then $\eta$ is regarded as a $\Cal V^2-$valued one from on the Hilbert manifold $\Cal V^1$. 
We immediately see that $\eta $ is a closed form. 
Hence the pull back $i^*_{\ss-style\E^1}\eta$ is also closed
where $i_{\ss-style\E^1}\: \E^1_s \to \Cal V^1$ denotes the inclusion. 
The distribution $Z$ is defined by the closed one form $i^*_{\ss-style\E^1}$, 
$$
Z_{\Phi } = \{ \, a \in \E^1_s \, |\, i^*_{\ss-style\E^1}\eta(a) =0\, \}.
\tag2-5$$
Hence we see that the distribution $Z$ is integrable. 
From the Frobenius theorem for Hilbert manifolds, there exists the integrable manifold $\CZ'$ such that 
$T_\Phi \CZ' = Z_\Phi$. We then see that $\CZ$ coincides with the integral manifold $\CZ'$. 
\qed\enddemo
Let $\Phi^0$ be a smooth closed element of $\E^1$. Since the orbit $\O$ is metrical, we have the
corresponding  smooth metric $g_{\ss-style\Phi^0}$. The vector bundle $E^i$ is defined from $\Phi^0$ as in
section one for each
$i=0,1,2$. Then we have the orthogonal projection $\pi_{\ss-style{E^i}}$, 
$$
\pi_{\ss-style{E^i}} \:  L^2_{s-1} ( X,\oplus_i\w^{p_i+i-1} )  \arrow L^2_{s-1} (E^i),
$$
for $i=0,1$.
We also denote by $d^*$ the adjoint operator with respect to the metric $g_{\ss-style\Phi^0}$,
We consider the complex:
$$
\CD
L^2_{s+1}(E^0) @>d_0>>L^2_s (E^1)@>d_1>>L^2_{s-1}(E^2).\qquad \quad (\#_{\ss-style\Phi^0})
\endCD
\tag2-6$$
Since $\O$ is elliptic, the complex $\#_{\ss-style\Phi^0}$ is an elliptic complex.
We then have the Laplacian  $\triangle_{\ss-style\#} =d_0 d_0^* + d_1^* d_1$ of the complex
$\#_{\ss-style\Phi^0}$ and the Hodge decomposition:
$$
L^2_s ( E^1)=\H^1(\#_{\ss-style\Phi^0})+ \triangle_{\ss-style\#} G_{\ss-style\#} L^2_s ( E^1),
\tag2-7$$where $\H^1(\#_{\ss-style\Phi^0})=
Ker \triangle_{\ss-style\#} $  
and $G_{\ss-style\#}$ denotes the Green operator.
We also see the the adjoint operator $d^*_0$(resp.  $d^*_1$) is given by $\pi_{\ss-style E^0}\circ d^*$
(resp. 
$\pi_{\ss-style E^1}\circ d^*$). \proclaim{Lemma 2-3}
Let $\Phi^0$ be a smooth closed element of $\E^1_s$. 
We define  a slice $S_{\ss-style\Phi^0}$ by 
$$
S_{\ss-style\Phi^0}=\{\, \Phi  \in \E^1_s \, |\, d\Phi =0, \,\pi_{\ss-style{E^0}} \circ d^* \Phi =0\, \} .
$$ 
Then $S_{\ss-style\Phi^0}\cap U$ is a Hilbert manifold with $T_{\Phi^0}S_{\ss-style\Phi^0}
=\H^1(\#_{\ss-style\Phi^0})$, where $U$ is a
sufficiently small neighborhood of $\E^1_s$ at $\Phi^0$.
\endproclaim
\demo{Proof} 
The slice $S_{\ss-style\Phi^0}$ is written as 
$$
S_{\ss-style\Phi^0}= \{\, \Phi \in \CZ\, |\, d_{{0}}  d_0^* \Phi =0\, \},
$$where $d_o^*$ denotes $ \pi_{E^0}\circ d^*$.
We define the map $F$ by 
$$\align
F\: &\CZ \arrow L^2_{s-2}(X, E^1),\\
& \Phi  \arrow  d_0  d_0^* \Phi
\endalign$$
The map $F$ is the map from the Hilbert manifold $\CZ$ to the closed sub space of $d_0$ exact forms 
$d_0L^2_{s-1}(E^0)$. 
The differential of $F$ at $\Phi^0$ is given by 
$$\align
dF_{\ss-style\Phi^0} \: T_{\Phi^0}&\CZ \arrow d_0L^2_{s-1}(E^0) \\
&a \arrow  d_0 d_0^* a,
\endalign$$
where $a \in T_{\Phi^0}\CZ= Z_{\Phi^0}$.
From the Hodge decomposition of the complex $\#_{\ss-style\Phi^0}$, 
we see that $dF_{\ss-style\Phi^0}$ is surjective. 
Hence from the implicit function theorem, the slice $S_{\ss-style\Phi^0}\cap U$ is 
a Hilbert manifold for a neighborhood $U$ at $\Phi^0$. 
The tangent space of $T_{\Phi^0}S_{\ss-style\Phi^0}$ is given by the kernel of the map
$dF_{\ss-style\Phi^0}$. From lemma 2-2 we see that the kernel of $dF_{\ss-style\Phi^0}$ is
$\H^1(\#_{\ss-style\Phi^0})$.
\qed\enddemo
\proclaim{Lemma 2-4}
We define a map $P_{_{S_{\Phi^0}}}$ by 
$$
P_{\ss-style{S_{\Phi^0}}}\: S_{\ss-style{\Phi^0}} \arrow \underset i\to \oplus
H^{p_i}_{\ss-style{dR}} (X),
$$
$$
P_{\ss-style{S_{\Phi^0}}}(\Phi ) = ( [\phi_1]_{\ss-style{dR}}, \cdots, [\phi_l]_{\ss-style{dR}} ),
$$
where $\Phi =(\phi_1,\cdots ,\phi_l)\in S_{\Phi^0}$ and $[\phi_i]_{\ss-style{dR}}$ is 
a class of
de Rham cohomology group represented by
$\phi_i$. Then $P_{S_{\Phi^0}}$ is injective on a small neighborhood at $\Phi^0$.
\endproclaim
\demo{Proof}
An open set of H$^1 (\#_{\Phi^0} )$ is a local coordinates of the slice $S_{\Phi^0}$. Hence the result
follows since the orbit $\O$ is topological.
\qed\enddemo
Let $C^1$Diff$_0$ be the identity component of $C^1-$diffeomorphisms from $X$ to $X$. 
We define Diff$^{s+1}_0$ by  
$$
\text{Diff}^{s+1}_0 = C^1\text{Diff}_0 \cap L^2_{s+1} (\text{ Diff }(X)).
$$
Then Diff$^{s+1}_0$ is the Hilbert Lie group and the action 
$$
\E^1_s \times \text{Diff}^{s+1}_0 \arrow \E^1_s
$$
is well defined (see [3] ).
We define $\M^s_{\ss-style{\O}} (X)$ by 
$$
\M^s_{\ss-style{\O}}(X)= \widetilde{\M}^s_{\ss-style{\O}} (X) /\text{Diff}_0^{s+1},
$$where
$$
\widetilde{\M}^s_{\ss-style{\O}} (X) = 
\{ \,\Phi \in \E^1_s \, |\, d\Phi =0 \, \}.
$$
Let $\pi$ be the natural projection 
$$
\pi\:\widetilde{\M}^s_{\ss-style{\O}} (X) \arrow \M^s_{\ss-style{\O}}(X).
$$
\proclaim{Lemma 2-5}
Let $S_{\Phi^0}$ be the slice through $\Phi^0$.
Then image $\pi ( S_{\Phi^0})$ is an open set in $\M^s_{\ss-style{\O}}(X)$.
\endproclaim
\demo{Proof} 
As in proof of lemma 2-3 it follows from implicit function theorem that a neighborhood of $\Phi^0$ is
homeomorphic to 
$S_{\Phi^0} \times V$, where $V$ is an neighborhood of Diff$_0$ at the identity. 
Hence we have the result.\qed
\enddemo
\proclaim{Lemma 2-6} 
Each element $\Phi$ of the slice $S_{\Phi^0}$ consists of smooth forms, i.e., 
$$\Phi = (\phi_1,\cdots,\phi_l), \quad\phi_i \in C^\infty (X,\w^{p_i}).
$$
So that is, 
$$
S_{\Phi^0} \subset \M_{\ss-style{\O}} (X). 
$$
\endproclaim
\demo{Proof}
The tangent space $T_{\Phi^0}S_{\Phi^0}$ is the Kernel of Laplacian $\triangle_{_\#}$. 
Hence from elliptic regularity we have 
$$
T_{\Phi^0}S_{\Phi^0}\subset C^\infty ( X, \underset i\to \oplus\w^{p_i}).
$$
Hence from implicit function theorem we see that 
$$
S_{\Phi^0} \subset C^\infty ( X, \underset i\to \oplus\w^{p_i}).
$$
\qed\enddemo
\proclaim{Proposition 2-7}
Let $\pi $ be the natural projection $\pi \:\til{\M}_{\ss-style{\O}}(X)\to \M_{\ss-style{\O}}(X)$.
We restrict the map $\pi$ to a slice $S_{\Phi^0}$. Then the restricted map $\pi|_{S_{\Phi^0}}$ to the image
$$
\pi|_{S_{\Phi^0}} \:S_{\Phi^0} \arrow \pi ( S_{\Phi^0})
$$
is a homeomorphism.
\endproclaim
\demo{Proof}
It is sufficient to show that $\pi|_{S_{\Phi^0}}$ is injective. 
We assume that $\pi (\Phi) = \pi (\Phi')$ for 
$\Phi,\Phi' \in S_{\Phi^0}$. 
It implies that there exists $f \in $Diff$_0$ such that 
$\Phi' = f^*\Phi$.
Since each class of de Rham cohomology group is invariant under the action of Diff$_0$, we have 
$$
[\Phi]_{\ss-style{dR}} = [\Phi']_{\ss-style{dR}}\in \underset i \to \oplus H^{p_i} (X).
$$
From lemma 2-4, the map $P_{S_{\Phi^0}}$ is injective for sufficiently small $S_{\Phi^0}$. 
Hence we have $\Phi =\Phi'$. 
\qed\enddemo
\proclaim{Proposition 2-8} 
The quotient $\M_{\ss-style{\O}}(X)$ is Hausdorff.
\endproclaim
\demo{Proof} 
We define $\til{\M}_s(X)$ by 
$$
\widetilde{\M}_s( X) = \{ \,\Phi \in \E^1_s\, |\, d\Phi =0 \, \}.
$$
Since $\O$ is metrical, we have the metric $g_\Phi$. Hence each tangent space 
$T_\Phi \widetilde{\M}(X)$ has the metric and 
$\til{\M}_s(X)$ is 
a Riemannian manifold. We also see that the action of Diff$_0$ on $\til{\M}_s(X)$ 
is isometric ( see [3] ).  Let $d$ be the distance of the Riemannian manifold 
$\til{\M}_s(X)$ and $\pi$  the natural projection
$\pi\: \widetilde{\M}_s(X)\to \M_{\ss-style{\O}}(X)$. Then we define $ d(
\pi(\Phi^1),\pi(\Phi^2) )$ by 
$$
d(\pi(\Phi^1) ,\pi(\Phi^2)) = \inf_{\ss-style{f,g \in \text{Diff}_0}} d ( f^*\Phi^1,
\,\,g^*\Phi^2 ),
$$
where $\Phi,\Phi' \in \til{\M}_s(X)$.
Since the action of $\Diff_0$ preserves the distance $d$,  
$$
d( \pi(\Phi^1),\pi(\Phi^2) ) = \inf_{\ss-style{f \in \Diff_0}}d ( f^*\Phi^1, \,\,\Phi^2).
$$
Hence we have triangle inequality , 
$$
\align
&\quad\quad  d(\pi (\Phi^1), \,\,\pi(\Phi^2) ) + d( \pi(\Phi^2),\,\,\pi(\Phi^3)) \\
=&\inf_{\ss-style{f \in \Diff_0}}d( f^*\Phi^1,\Phi^2) + 
\inf_{\ss-style{g\in \Diff_0} }d(\Phi^2,
g^*\Phi^3)\\ 
\leqq &\inf_{\ss-style{f,g\in \Diff_0}} d( f^*\Phi^1,\,\,g^*\Phi^3) = d(
\pi(\Phi^1),\,\pi(\Phi^3)).
\endalign
$$
We shall that $d$ is a distance of $\M_{\ss-style{\O}}(X)$. 
We assume that the distance $d(\,\pi(\Phi^0),\,\,\pi(\Phi)\,)=0$. 
Then ${\inf_{\ss-style{f\in\Diff_0}}d(\Phi^0, f^*\Phi)=0}$. Hence $f^*\Phi$ is in a small
neighborhood 
$U$ at $\Phi^0$.  As in lemma 2-1
  a neighborhood $U$ of $\wtil{\M}(X)$ at $\Phi^0$ is homeomorphic to a product $
V\times S_{\Phi^0}$, where $V$ is a neighborhood of \Diff$_0$ at the identity. 
We define a distance $d_{\Phi^0}$ on the cohomology group  $\oplus_i H^{p_i}_{\ss-style{dR}}(X)$ by using 
the harmonic representation with respect to the metric
$g_{\Phi^0}$.  From lemma 2-4, we have the injective map 
$$
P\: S_{\Phi^0}\to \oplus_i H^{p_i}_{\ss-style{dR}}(X).
$$
Then we see that 
$$
\inf\Sb {\ss-style{f\in\Diff_0,}}\\{ \ss-style{f^*\Phi\in U}}\endSb d(\Phi^0, f^*\Phi)=
 C\,\, d_{\ss-style{\Phi^0}}(
P(\Phi^0),\,\,P(\Phi)),
$$
where $C$ is a constant. (Note that  the action of \Diff$_0$ preserves a class of $\oplus_i
H^{p_i}_{\ss-style{dR}}(X)$.) Hence from our assumption we have $P(\Phi) = P(\Phi^0)$. 
since $P$ is injective, $\pi(\Phi)=\pi(\Phi^0)$. Hence $d$ is a distance on
$\M_{\ss-style{\O}}(X)$.  
\qed\enddemo
\demo{Proof of theorem 1-6}
Since $\O$ is elliptic, each slice $S_{\Phi}$ is smooth from lemma 2-3. 
Each slice $S_{\Phi}$ is  coordinates of $\M_{\ss-style{\O}}(X)$ since $\O$ is topological
from lemma 2-5.  
Each slice is homeomorphic to an open set of the cohomology group H$^1 (\#_\Phi )$. 
Since $\O$ is metrical, $\M_{\ss-style{\O}}(X)$ is Hausdorff and each slice is canonically constructed.
\qed\enddemo
\demo{Proof of theorem 1-7}
Each slice is a local coordinates of $\M_{\ss-style{\O}}(X)$. Then the result follows from lemma 2-4.
\qed\enddemo
Since $\O$ is metrical, we have the metric $g_{\Phi}$ for each $\Phi \in \E^1$.
Hence the metric $g_{\Phi}$ defines the metric on each tangent space $E^1 =
T_\Phi \E^1$.  So $\E^1$ can be considered as a Riemannian manifold.
Then we see that the action of Diff$_0$ on  $\E^1$ is isometry. 
Let $I({\Phi})$ be the isotropy group of Diff$_0$ at $\Phi$, 
$$
I({\Phi}) = \{\, f \in \text{Diff}_0(X) \, |\, f^*\Phi =\Phi\, \}.
$$ 
Let $S_{\Phi^0}$ be a slice at $\Phi^0$. Then we shall compare $I_{\Phi^0}$ to 
other isotropy group $I_{\Phi}$ for $\Phi \in S_{\Phi^0}$.
\proclaim{Theorem 1-8}
Let $I({\Phi}^0)$ be the isotropy group of Diff$_0(X)$ at $\Phi^0$ and $S_{\Phi^0}$ the slice at
$\Phi^0$.  Then $I_{\Phi^0}$ is a  subgroup of the isotropy group $I_{\Phi}$ for each $\Phi\in
S_{\Phi^0}$. ( We take $S_{\Phi^0}$ sufficiently small for necessary.)
\endproclaim
\demo{Proof of theorem 1-8}
From definition of $S_{\Phi^0}$, the slice $S_{\Phi^0}$ is invariant under the action of 
$I_{\Phi^0}$. The map $P|_{S_{\Phi^0}} \: S_{\Phi ^0}\to \underset i \to \oplus $H$^{p_i} (X)$ is 
locally injective. Since the action of Diff$_0$ preserves each class of de Rham cohomology group, we can take a
sufficiently small 
$S_{\Phi^0}$ such that the action of $I_{\Phi^0}$ is trivial on the slice $S_{\Phi^0}$.
Hence $I_{\Phi^0}$ is a subgroup of the isotropy group $I_{\Phi}$ for each 
$\Phi \in S_{\Phi^0}$.
\qed\enddemo
\demo{Proof of theorem 1-9}
The slice $S_{\Phi^0}$ is local coordinates of $\M_{\ss-style{KE}}(X)$ and 
invariant under the action of $\Diff(X)$. 
Hence the moduli space $\wtil{\M}_{\ss-style{KE}}(X)/\Diff(X)$ is locally homeomorphic to 
the quotient space $S_{\Phi^0}/I$, where $I$ is the isotropy, 
$$
I = \{ \, f \in \Diff(X)\,| f^* \Phi^0 =\Phi^0 \, \}.
$$
As in proof of proposition 2-8, $\Diff(X)$ acts on $S_{\Phi^0}$ isometrically. 
Hence we see that there is an open set $V$ of $T_{\ss-style{\Phi^0}} S_{\ss-style{\Phi^0}}$
with the action of
$I$  such that the quotient $V/I$ is homeomorphic to $S_{\Phi^0}/I$. 
Since $T_{\Phi^0}S_{\Phi^0}$ is isomorphic to $H^1(\#_{\ss-style {\Phi^0}})$ and the action of
$I$ on $H^1(\#_{\ss-style {\Phi^0}})$ is a isometry with respect to $g_{\ss-style {\Phi^0}}$. The
action of $I$ preserves integral  cohomology class. Hence from lemma 2-4 we see that  $V/I$ is 
the quotient by a finite group. 
\qed\enddemo
\head 
\S3-1. Calabi-Yau structures 
\endhead
Let $V$ be a real vector space of $2n$ dimensional. 
We consider the complex vector space $V\otimes \C$ and 
a complex form $\Ome \in \w^n V^*\otimes\C$. 
The vector space ker\,$\Ome$ is defined as 
$$
Ker\, \Ome = \{ \, v \in V\otimes \C \, |\, i_v \Ome =0 \, \},
$$
where $i_v$ denotes the interior product.
\proclaim{Definition 3-1-1 (Calabi-Yau structures) } 
A complex $n$ form $\Ome$ is a Calabi-Yau structure on $V$ if 
$\dim_\C $Ker\, $\Ome =n$ and Ker\, $\Ome \cap \ol{Ker\,\Ome}= \{0\}$,
where $\ol{Ker\, \Ome}$ is the conjugate vector space.
\endproclaim
We denote by $\A_{\ss-style{CY}}(V)$ the set of Calabi-Yau structures on $V$. 
We define the almost complex structure $I_\Ome$ on $V$ by 
$$
I_\Ome (v) = 
\cases
\sqrt{-1}v&\quad \text{ if  } v\in Ker\, \Ome,\\
-\sqrt{-1}v&\quad \text{ if }v \in \ol{Ker\, \Ome}.
\endcases
$$
So that is, Ker\, $\Ome= T^{1,0}V$ and $\ol{Ker \Ome}=T^{0,1}V$ and
$\Ome$ is a non-zero $(n,0)$ form on $V$ with respect to $I_\Ome$.
Let $ \Cal J$ be the set of almost complex structures on $V$. Then 
$\A_{\ss-style{CY}}(V)$ is the $\C^*-$bundle over$\Cal J$. 
We denote by $\rho$ 
the action of the real general linear group $G=GL (V)$ on the complex $n$ forms, 
$$
\rho \: \text{GL}(V) \arrow \text{End}\,(\w^n (V\otimes\C)^* ).
$$
Since $G$ is a real group, $\A_{\ss-style{CY}}(V)$ is invariant under the action of $G$. 
Then we see that the action of $G$ on $\A_{\ss-style{CY}}(V)$ is transitive. 
The isotropy group  $H$ is defined as 
$$
H =\{\, g \in G\, |\, \rho_g \Ome = \Ome \, \}.
$$
Then we see $H=$SL$(n, \C)$. Hence the set of Calabi-Yau structures $\A_{CY
}(V)$ is the homogeneous space, 
$$
\A_{\ss-style{CY}}(V) = G/H = GL(2n, \R) /SL(n,\C).
$$
(Note that the set of almost complex structures $\Cal J=$GL$(2n,\R)/$GL$(n,\C)$. )
An almost complex structure $I$ defines a complex subspace $T^{1,0}$ of dimension
$n$. Hence we have the map $\Cal J \arrow $Gr$(n, C^{2n})$. We also have the map 
from $\A_{\ss-style{CY}}(V)$ to the tautological line bundle $L$ over the Grassmannian 
Gr$(n,\C^{2n})$ removed $0-$section.
Then we have the diagram: 
$$
\CD 
\A_{\ss-style{CY}}(V) @>>> L\backslash{0}\\
@VV\C^*V @VVV\\
\Cal J @>>>Gr(n,\C^{2n} )
\endCD
$$
$\A_{\ss-style{CY}}(V)$ is embedded as a smooth submanifold in $n-$forms $\w^n$. 
This is Pl\"ucker embedding described as follows, 
$$
\CD
\A_{\ss-style{CY}}(V) @>>> L\backslash{0}@>>>\w^n\backslash\{0\}\\
@VV\C^*V @VVV@VVV\\
\Cal J @>>>Gr(n,\C^{2n} )@>>>\CP^n.
\endCD
$$
Hence the orbit $\O_{\ss-style{CY}}=\A_{\ss-style{CY}}(V)$ is a submanifold in $\w^n$ 
defined by Pl\"uker relations.
Let $X$ be a real compact manifold of $n$ dimensional. 
Then we have the $G/H$ bundle $\A_{\ss-style{CY}}(X)$ over $X$ as in section 1. 
We denote by $\E=\E^1_{\ss-style{CY}}$ the set of smooth global section of $\A_{\ss-style{CY}}(X)$. 
Then we have the almost complex structure $I_\Ome$ corresponding to 
$\Ome\in \E^1$. Then we have 
\proclaim{Lemma 3-1-2}
If $\Ome \in \E^1$ is closed, then the almost complex structure $I_\Ome$ is 
integrable.
\endproclaim
\demo{Proof}
Let $\{\theta_i\}_{i=1}^n$ be a local basis of $\Gam(\w^{1,0})$ with respect to 
$\Ome$. From Newlander-Nirenberg's theorem it is sufficient to show that 
$d\theta_i \in \Gam (\w^{2,0}\oplus \w^{1,1} )$ for each $\theta_i$. 
Since $\Ome$ is of type $\w^{n,0}$, 
$$
\theta_i \w \Ome =0.
$$
Since $d\Ome=0$, we have 
$$
d\theta_i \w \Ome =0.
$$
Hence $d\theta_i \in \Gam(\w^{2,0}\oplus \w^{1,1} )$.
\qed\enddemo
Then we have the moduli space of Calabi-Yau structures on $X$, 
$$
\M_{\ss-style{CY}}(X) = \{\, \Ome \in \E^1_{\ss-style{CY}}\, |\, d\Ome =0 \, \}/
\text{Diff}_0(X).
$$ 
From lemma 3-1-2 we see that $\M_{\ss-style{CY}}(X)$ is the $\C^*-$bundle over 
the moduli space of integrable complex structures on $X$ with trivial canonical
line bundles.
\proclaim{Proposition 3-1-3}
The orbit $\O_{\ss-style{CY}}$ is elliptic. 
\endproclaim
\demo{Proof}
Let $\w^{p.q}$ be $(p,q)-$forms on $V$ with respect to  a closed form $\Ome^0 \in
\A_{\ss-style{CY}}(V)$.  In this case we see that 
$$\align 
E^0 &= \w^{n-1,0}\\
E^1 &=\w^{n,0}\oplus \w^{n-1,1}\\
E^2 &=\w^{n,1}\oplus \w^{n-1,2}.
\endalign
$$
Hence we have the complex :
$$\CD
 \w^{n-1,0} @>\w u >>\w^{n,0}\oplus \w^{n-1,1}@>\w u >> \w^{n,1}\oplus\w^{n-1,2},
\endCD
$$
for $ u \in V$. 
Let $ a=(x, y )$ be an element of $E^1$. We assume that $u\w a=0$. 
Then since de Rham complex is elliptic, there is an element $b =(z,w)\in
\w^{n-1,0}\oplus\w^{n-2,1}$ such that $a = u\w b$. Hence we have 
$$\align
x& = u^{1,0}\w z,\tag1\\
y &= u^{0,1}\w z + u^{1,0}\w w, \tag2\\
0&= u^{0,1}\w w \tag 3
\endalign $$
From the equation (3) and ellipticity of Dolbeault complex, there is an element
$\gam
\in \w^{n-2,0}$ such that $ w = u^{0,1} \w \gam$.  Hence $y = u^{0,1}\w z + u^{1,0}\w
(u^{0,1}\w \gam )= u ^{0,1} \w \h{z}$, where $\h{z} =z - u^{1,0}\w \gam\in \w^{n-1,0}$. 
Hence $ a = (u^{1,0}\w \h{z}, u^{0,1}\w \h{z} )$ for $\h{z} \i \w^{n-1,0}$. 
\qed\enddemo
\proclaim{Proposition 3-1-4}
Let $I_\Ome$ be the complex structure corresponding to $\Ome \in \E^1$. 
If $\partial \ol{\partial}$ lemma holds for the complex manifold $(X, I_\Ome)$,
$\Ome $ is topological. In particular, $(X, I_\Ome)$ is K\"ahkerian, $\Ome$ is 
topological.
\endproclaim
\demo{Proof}
As in proof of proposition 3-1-2 the complex \#$_\Ome$ is given as
$$
\CD 
\Gam( \w^{n-1,0} ) @>d>>\Gam ( \w^{n,0}\oplus\w^{n-1,1} ) @>d>>\Gam (\w^{n,1}
\oplus\w^{n-1,2} ).
\endCD
$$
Then we have the following double complex:

$$\CD 
\Gam(\w^{n,0} )@>\ol{\pa}>>\Gam(\w^{n,1})@>\ol{\pa}>>\\
@A\pa AA @ A\pa AA @.\\
\Gam(\w^{n-1,0})@>\ol{\pa}>>\Gam(\w^{n-1,1})@>\ol{\pa}>>\\
@A\pa AA @A\pa AA @. \\
\Gam(\w^{n-2,0})@>\ol{\pa}>>\Gam(\w^{n-2,1})@>\ol{\pa}>> 
\endCD$$
Let $a = (x,y)$ be a closed element of $\Gam(\w^{n,0})\oplus\Gam(\w^{n-1,1})$. 
We assume that $a = db$. Then $b=(z,w) \in \Gam( \w^{n-1,0})\oplus\Gam (\w^{n-2,1} )$,
satisfying 
$$
\align 
x &= \pa z \\
y&= \ol{\pa}z +\pa w \\
0&= \ol{\pa}w.
\endalign
$$
Since $\pa w \in \Gam (\w^{n-1,1})$ is $\ol{\pa}$ closed, 
We apply $\pa \ol{\pa}$ lemma to $\pa w$. 
Then there is an element $\gam \in \Gam (\w^{n-2,0})$ such that 
$$
\pa w = \pa\ol{\pa}  \gam = - \ol{\pa} \pa\gam.
$$
Hence $a$ is written as 
$$
a = d ( z - \pa \gam).
$$
It implies that the map $p\: H^1 (\#) \arrow H^n (X)$ is injective.
\qed\enddemo
\demo{Remark} 
We define $F^{n-1}\w^* $ 
$$
F^{n-1}\w^m = 
\bigoplus\Sb p+q =m \\
 p \geq n-1\endSb
 \Gam ( \w^{p,q} ).
$$
Then the complex \#$_\Ome$ is 
$$
\CD
F^{n-1} \w^{n-1} @>d>>F^{n-1}\w^n @>d>>F^{n-1}\w^{n+1}.
\endCD$$
Hence $\Ome$ is topological if and only if we have the following: 
$$
F^{n-1}H^n (X) = H^n (F^{n-1} \w^*).
$$
\enddemo
Hence from proposition 3-1-3, we have the smooth slice $S_\Ome$ corresponding to each
$\Ome$ and $S_\Ome$ is the space of deformations of $\Ome$. 
However $\O_{\ss-style{CY}}$ is not metrical and
the moduli space
$\M_{\ss-style{CY}}(X)$ is not Hausdorff in general. 
In fact, It is known that $\M_{\ss-style{CY}}(X)$ is not Hausdorff for a $K3$ surface.
Hence in order to obtain a Hausdorff moduli space, we must introduce 
extra geometric structures. The most natural structure is K\"ahler-Einstein structure 
on a Calabi-Yau manifold.
\head \S3-2. K\"ahler-Einstein structures
\endhead
Let $V$ be a real vector space of $2n$ dimensional. 
We consider a pair $\Phi=( \Ome, \ome)$ of a Calabi-Yau structure $\Ome$ and 
a real symplectic structure $\ome$ on $V$, 
$$\align
\Ome &\in \A_{\ss-style{CY}}(V), \\
 \ome &\in \w^2 V^*, \quad\overset n \to{\overbrace 
{\ome\w
\cdots \w \ome} }\neq 0.
\endalign$$
We define $g_{\Ome,\ome}$ by 
$$
g_{\Ome,\ome}(u,v)=\ome(I_\Ome u, v), 
$$
for $u,v \in V$.
\proclaim{Definition 3-2-1(K\"ahler-Einstein structures )} 
A K\"ahler-Einstein structure on $V$ is a pair $\Phi=( \Ome,\ome) $ such that 

$$\align
&\Ome \w \ome =0 , \quad \ol{\Ome}\w \ome =0
\tag 1\\
&\Ome \w \ol{\Ome} = 
c_n \,\overset n\to{\overbrace{\ome\w \cdots \w\ome} 
}
\tag 2\\
&g_{\Ome,\ome}\text{ is positive definite.}\tag3
\endalign$$
where $c_n$ is a constant depending only on $n$,.i.e, 
$$
c_n =(-1)^{\frac{n(n-1)}2}\frac{2^n}{i^n n!}.
$$
\endproclaim
From the equation (1) we see that $\ome$ is of type $\w^{1,1}$ with respect to 
the almost complex structure $I_\Ome$. The equation (2) is called
Monge-Amp$\grave{e}$re equation. 
\proclaim{Lemma 3-2-2}
Let $\A_{\ss-style{KE}}(V)$ be the set of K\"ahler-Einstein structures on $V$. 
Then There is the  transitive action of $G=$GL$( 2n, \R)$ on $\A_{\ss-style{KE}}(V)$ 
and $\A_{\ss-style{KE}}(V)$ is the homogeneous space 
$$
\A_{\ss-style{KE}}(V)= GL(2n,\R) / SU(n).
$$ 
\endproclaim
\demo{Proof}Let $g_{\Ome,\ome}$ be the K\"ahler metric. Then we have a unitary
basis  of $TX$. Then the result follows from (1),(2).
\qed\enddemo
Hence the set of K\"ahler-Einstein structures on $V$ is the orbit $\O_{\ss-style{KE}}$, 
$$
\O_{\ss-style{KE}} \subset \w^n (V\otimes \C)^*\oplus \w^2 V^*.$$
\proclaim{Theorem 3-2-3} 
The orbit $\O_{\ss-style{KE}}$ is metrical, elliptic and topological.
\endproclaim
\demo{Proof}
From lemma 3-2-2 the isotropy group is SU$(n)$. Hence $\O_{\ss-style{KE}}$ is metrical. 
At first we shall show that $\O_{\ss-style{KE}}$ is elliptic. 
Let $(\Ome^0,\ome^0)$ be an element of $\A_{\ss-style{KE}}(V)$. 
Then we have the vector space $E^0(V)=E^0_{\ss-style{KE}}(V)$ by 
$$
E^0_{\ss-style{KE}}(V) = \{\, (i_v \Ome^0, i_v\ome^0 ) \in \w^{n-1}_\C\oplus \w^{n-1}\,|
\, v \in V\, \}
$$
The vector space $E^1(V)=E^1_{\ss-style{KE}}(V)$ is the set of 
$ (\a,\b) \in \w^n_\C \oplus \w^2$
satisfying equations 
$$\align
&\a \w \ome^0 + \Ome^0 \w \b =0 , \\
&\a \w \ol{\Ome^0} + \Ome^0 \w \ol{\a} = n c_n \b \w (\ome^0)^{n-1} 
\tag 4
\endalign
$$
We assume that $u\w \a =0 , u\w \b =0 $ for some non zero vector $u \in V$. 
Then since the orbit $\O_{\ss-style{CY}}$ is elliptic, 
$(\a ,\b)$ is given as 
$$
\a = u \w s ,\quad \b = u\w t,
\tag5$$
form some $s \in \w^{n-1,0}_{I_\Ome^0}$ and $t \in \w^1$. 
Hence $s,t$ are written as 
$$
s = i_{v_1}\Ome^0,\quad t =i_{v_2}\ome^0,
\tag6$$
for some $v_1, v_2 \in V$. 
Set $v=v_1 -v_2$. Then from (4),(5) and (6) using (1),(2),
we have 
$$
\align 
&u \w ( i_v \ome^0) \w \Ome^0 =0 
\tag 7
\\
&u\w ( i_v \ome^0) \w (\ome^0)^{n-1} =0.
\tag 8
\endalign 
$$
Form (7) we have 
$$
u\w i_v \ome^0 \in \w^{2,0} \oplus \w^{1,1}.
\tag 9$$
We also have from (8)
$$
u \w i_v \ome^0 \in \w^{2,0}\oplus \w^{0,2}.
\tag 10$$
Since $u\w i_v \ome^0 $ is a real form, we
see from (9),(10) that 
$$
u\w i_v \ome^0 =0.
\tag 11
$$
Hence $(\a,\b)$ is given as 
$$\align 
&\a = u \w i_{v_1}\Ome^0,\\
&\b=u\w i_{v_2}\ome^0 = u\w i_{v_1 }\ome^0 -u\w i_v \ome^0 =u\w
i_{v_1}\ome^0.
\tag 12
\endalign$$
From (12) we see that 
the complex 
$$\CD
E_{\ss-style{KE}}^0 (V) @>\w u >> E_{\ss-style{KE}}^1(V) @> \w u >> E_{\ss-style{KE}} ^2(V)
\endCD$$
is elliptic.
Next we shall show that $\O_{\ss-style{KE}}$ is topological. 
Let $(\a,\b )$ be an element of $\Gam( E_{\ss-style{KE}}^1)$. 
We assume that $\a$ and $\b$ are exact forms respectively.
Then since $\O_{\ss-style{CY}}$ is topological, we have 
$ \a = ds,\b =dt$ for some 
$s \in  \Gam( \w^{n-1}_\C),$ $ t \in \Gam(\w^2)$. 
Hence  $s,t$ are written as 
$$
s=i_{v_1}\Ome^0, t =i_{v_2}\ome^0,
\tag13$$
for some $v_1,v_2 \in \Gam (TX)$.
Then from equations (4),(13) using (1),(2) we have 
$$
\align 
&d(i_v \ome^0) \w \Ome^0 =0 \tag 14\\
&d(i_{v}\ome^0)\w (\ome^0)^{n-1} =0,
\tag 15\endalign
$$
where $v = v_1 -v_2$.
From (14),(15) we have 
$$
\align
&d i_v \ome^0 \in \Gam ( \w^{2,0} \oplus \w^{1,1} )\tag 16\\ 
&d_v \ome^0 \in \Gam( \w^{2,0}\oplus \w^{0,2} ).
\tag 17 
\endalign
$$
Since $di_v \ome^0$ is real, we have from (16), (17) 
$$
di_v \ome^0 =0.
\tag 18$$
Hence $\b= di_{v_2}\ome^0 = d i_{v_1}\ome^0 -d i_v \ome^0$. 
Then from (18) we see that 
$$
\a = d i_{v_1} \Ome^0 , \quad \b =di_{v_1} \ome^0.
$$
Hence the map $p \: H^1( \#) \arrow H^n (X, \C) \oplus H^2 (X.\R)$ 
is injective.
\qed\enddemo
Hence from theorem 1-6 in section 1 we have the following:
\proclaim{Theorem 3-2-4}
Let $\M_{\ss-style{KE}}(X)$ be the moduli space of K\"ahler-Einstein structures over $X$, 
$$
\M_{\ss-style{KE}} (X) = \til{\M}_{\ss-style{KE}}(X)/\text{\rm Diff}_0(X),
$$
where 
$$
\til{\M}_{\ss-style{KE}}(X) =
\{ \, ( \Ome,\ome) \in \E^1_{\ss-style{KE}}\, |
\, d \Ome =0 , d\ome =0 \, \} .
$$
Then $\M_{\ss-style{KE}}(X)$ is a smooth manifold. 
Let $\pi $ be the natural projection 
$$\pi\: \til{\M}_{\ss-style{KE}}(X) \arrow 
\M_{\ss-style{KE}}(X).$$
Then coordinates of $\M_{\ss-style{KE}}(X)$ at each $(\Ome,\ome) \in \til{\M}_{\ss-style{KE}}(X)$ 
is canonically given by an open ball of the cohomology group $H^1 (\#)$. 
\endproclaim
We have the Dolbeault cohomology group $H^{p,q} (X)$ with respect to each $\Ome
$. Then we have 
\proclaim{Theorem 3-2-5} 
The cohomology group $H^1(\#)$ is the subspace of 
$H^n(X,\C) \oplus H^2 (X,\R)$ which is defined by equations 
$$\align
&\a \w \ome+ \Ome \w \b =0 , \\
&\a \w \ol{\Ome} + \Ome \w \ol{\a} = n c_n \b \w\ome^{n-1}, 
\tag 19
\endalign
$$
where $\a \in H^n (X,\C) , \b \in H^2 (X,\R)$.
\endproclaim
Let $P^{p,q}(X)$ be the primitive cohomology group with respect to
$\ome$. Then we have Lefschetz decomposition, 
$$
\a = \a^{\ss-style{n,0}} + \a^{\ss-style{n-1,1}} +\a ^{\ss-style{n-2,0}}\neg\w
\ome 
\in P^{n,0}(X) \oplus P^{n-1,0}(X) \oplus P^{n-2,0}(X)\negthinspace\w \ome .
$$
$$
\b = \b^{\ss-style{2,0}} + \b^{\ss-style{1,1}}+ \b^{\ss-style{0,0}}\neg\w \ome +
\b^{\ss-style{0,2}} \in
 P^{2,0}(X)\oplus P^{1,1}(X) \oplus P^{0,0}(X)\negthinspace\w \ome \oplus 
P^{0,2}(X).
$$
Then equation (19) is written as 
$$\align 
&\a^{\ss-style{n-2,0}}\w \ome \w \ome +\Ome \w \b^{\ss-style{0,2}} =0, \\
&\a^{\ss-style{n,0}}\w \ol{\Ome} = nc_n \b^{\ss-style{0,0}}\ome^n
\endalign$$
We see that
$\a^{\ss-style{n,0}} \in P^{n,0}(X)$ and $\b^{\ss-style{0,0}}\in P^{0,0}(X)$ are
corresponding to  the deformation in terms of
constant multiplication: 
$$
\Ome \arrow t \Ome, \quad
\ome \arrow s \ome
$$
If a K\"ahler class $[\ome]$ is not invariant under 
a deformation, such a deformation corresponds to an element of 
$\b^{2,0}$ and $\a^{n-2,0}$.  This is in the case of Calabi family of
hyperK\"ahler manifolds, i.e.,  Twistor space gives such a deformation.
It must be noted that there is no relation between $\a^{n-1,1}\in P^{n-1,1}(X)$ and 
$\b^{1,1}(X)\in P^{1,1}(X)$.
We have from theorem 1-7 in section 1,
\proclaim {Theorem 3-2-6}
The map $P$ is locally injective, 
$$
P\: \M_{\ss-style{KE}}(X) \arrow H^n (X,\C) \oplus H^2(X,\R).
$$
\endproclaim
We also have from theorem 1-8 in section 1,
\proclaim{Theorem 3-2-7} 
Let $I(\Ome,\ome)$ be the isotropy group of $(\Ome,\ome)$, 
$$
I(\Ome,\ome) = \{ \, f \in \text{\rm Diff}_0(X)\, |\, f^*\Ome = \Ome , \, f^*
\ome =
\ome \, \}.
$$
We consider the slice $S_0$ at $\Phi^0=( \Ome^0, \ome^0)$. 
Then the isotropy group $I( \Ome^0,\ome^0)$ is a subgroup of 
$I(\Ome,\ome)$ for each $( \Ome,\ome) \in S_0$.
\endproclaim

We define the map $P_{H^2}$ by 
$$
P_{H^2}\: \M_{\ss-style{KE}}(X) \arrow \Bbb P(H^2 (X)),
$$where
$$
P_{H^2} ( [\Ome,\ome] ) \arrow [\ome]_{dR} \in \Bbb P (H^2 (X)),
$$
$\Bbb P(H^2(X))$ denoted the projective space $(H^2(X)-\{0\})/ \R^*$. 
Then we have 
\proclaim{Theorem 3-2-8} 
The inverse image $P^{-1}_{H^2}([\ome]_{dR})$ 
is a smooth manifold. 
\endproclaim
\demo{Proof}
From theorem 3-2-5 and theorem 3-2-6 the differential of the 
map $P_{H^2}$ is surjective. Hence from the implicit function theorem 
$P^{-1}_{H^2}([\ome]_{\ss-style{dR}})$ is a smooth manifold.
\enddemo
\demo{Remark} 
$P^{-1}_{H^2}([\ome]_{\ss-style{dR}})$ is 
the $\C^*$ bundle over the moduli space of polarized Calabi-Yau manifolds [4].
\enddemo
\head\S4.  Special lagrangians and K\"ahler-Einstein structures
\endhead
\subhead \S4-1 \endsubhead
Let $\Phi^0=(\Ome^0,\ome^0)$ be a K\"ahler-Einstein structure  on
a compact Calabi-Yau manifold $X$ and 
$i_{\ss-style M}\: M\subset X$ a special lagrangian submanifold of $X$ with
respect to 
$(\Ome^0,\ome^0)$, .i.e, 
$$
i_{_M}^* (\Ome^0)^{Im}=0,\quad i_{_M}^* \ome^0=0,
$$
where $(\Ome)^{Im}$ is the imaginary part of the complex form $\Ome^0$. 
We assume that
$M$ is a compact $n$ dimensional real manifold.
We denote by Diff$_0(X, M)$ the subgroup of Diff$_0 (X)$ 
preserving the submanifold $M$,
$$
\text{Diff}_0(X,M)=\{ \, f \in \text{\rm Diff}_0 (X)\, |\, f (M) =M\, \},
$$
where Diff$_0(X)$ is the identity component of the group of diffeomorphisms
of $X$. 
\proclaim{Definition 4-1-1}
We define the relative moduli space
$\M_{_{\neg K\neg E}}\neg(X,M)$ by 
$$
\M_{_{\neg K\neg E}}\neg(X,M)=
\wtil{\M}_{_{\neg K\neg
E}}\neg(X,M)/\text{\Diff}_0(X,M),
$$
where $\wtil{\M}_{_{\neg K\neg E}}\neg(X,M)$ is given as
$$
\wtil{\M}_{_{\neg K\neg E}}\neg(X,M)=
\{\, \Phi=( \Ome,\ome)\in \E^1_{_{\neg K\neg E}}\neg (X)\, |
\,d\Phi=0,\,i_{_M}^*\Ome^{Im}=0, \, i_{_M}^*\ome =0 \, \}.
$$
\endproclaim
For simplicity we use the following notation:
\proclaim{Notation 4-1-2}
Let $\Phi=(\Ome,\ome)$ be a K\"ahler -Einstein structure. Then a pair 
$(\phi,\psi)$ 
denotes   
$\Phi$, where 
$$
\phi =\Ome^{Re} ,\quad \psi =(\Ome^{Im},\ome)
$$
where $\Ome^{Re}$ is the real part of $\Ome$.
\endproclaim
Then $\M_{_{\neg K\neg E}}\neg(X,M)=\M(X,M)$ is rewritten as 
$$
\M_{_{\neg K\neg E}}\neg(X,M) = \{ \Phi=(\phi,\psi)\in \E^1_{_{\neg K\neg E}}\neg(X) \, |
\, d\Phi=0, i_{_M}^*\psi =0\, \}/\text{Diff}_0(X,M)
$$
We consider a special lagrangian $M'$ with respect to a K\"ahler-Einstein 
structure $\Phi=(\Ome,\ome)$ on $X$. Then there is the action of
Diff$_0(X)$  on the set of pairs $(M',\Phi)$ by 
$$
f (M',\Phi)  =( f^{-1}(M),f^*\Phi),
$$
for $f \in $Diff$_0$. Hence we have the moduli space of pairs of 
special lagrangian submanifolds and K\"ahler-Einstein structures on $X$:
\proclaim{Definition 4-1-3} 
$$
\Cal P = \{\, ( \Phi, M) \, |\, \Phi \in \wtil{\M}_{_{\neg K\neg E}}\neg(X), 
M\subset X, \text{ a special lagrangian } \, \}/ \text{\rm Diff}_0(X),
$$
where 
$$\wtil{\M}_{_{\neg K\neg E}}\neg(X)=\{\, \Phi=( \Ome,\ome) \in \E^1_{_{\neg K\neg E}}\neg(X)\, |\, d\Phi =0
\,\}.$$
\endproclaim
\proclaim{Lemma 4-1-4}
$\M_{_{\neg K\neg E}}\neg(X,M)$ is a connected component of the
moduli space 
$\Cal P$.
\endproclaim
\demo{Proof} 
There is the natural map $\wtil{\M}_{_{\neg K\neg E}}\neg(X,M) \to \Cal P$.
From Definition 4-1-1 and 2, we see this map is injective. 
Let $M'$ be  a special lagrangian submanifold with respect
to  a K\"ahler-Einstein structure $\Phi'$. If the class $[\Phi',M'] \in \Cal
P$ belongs to the same connected component as to the one of  the class
$[M,\Phi]$, we have 
$$
M' = f(M),
$$
for some $f \in \Diff_0(X)$. Hence the map $\wtil{\M}_{_{\neg K\neg
E}}\neg(X,M) \to
\Cal P$ is bijective map to the connected component of $\Cal P$. 
\qed\enddemo
As in section 3,
we have vector bundles $E^0_{\neg_X}, E^1_{\neg_X}$ and
$E^2_{\neg_X}$ for each K\"ahler-Einstein structure $\Phi^0
=(\Ome^0,\ome^0)$ by  
$$
\align
&\Gam_{\neg_X}\neg(E^0_{\neg_X})= \{\, (i_v
\Ome^0, \,i_v \ome^0 ) \in 
\Gam_{\neg_X}\neg(\, \w^{n-1}_\C\oplus\w^1) 
\,|\, v \in \Gam_{\neg_X}\neg( TX)\, \}\\
&\Gam_{\neg_X}\neg(E^1_{\neg_X})= \{ \,(\theta\w i_v
\Ome^0,\theta\w i_v\ome^0 )\in 
\Gam_{\neg_X}\neg( \w^n\C\oplus\w^2 )\,|\,\theta\in
\Gam_{\neg_X}(\w^1),  v \in
\Gam_{\neg_X}\neg(TX)\,\} \\
&\Gam_{\neg_X}\neg(E^2_{\neg_X}) =\{\, (\theta\w\a,\theta\w\b )\in
\Gam_{\neg_X}\neg( \w_\C^{n+1}\oplus\w^3) \,|\, \theta \in
\Gam_{\neg_X}\neg(\w^1),(\a,\b)
\in
\Gam_{\neg_X}\neg(E^1_{\neg_X})\, \},
\endalign
$$we have the complex $\#_{_X}$= \#$_{_{\Phi^0}}$, 
$$
\CD
\Gam_{\neg_X}\neg(E^0_{\ss-style X})
@>d_{\ss-style{X,0}}>>\Gam_{\neg_X}\neg(E^1_{\ss-style X})
@>d_{X,1}>>
\Gam_{\neg_X}\neg(E^2_{\ss-style X})
\endCD
\tag \#$_{_{\Phi^0}}$ 
$$
By using Notation $\Phi^0=(\phi^0,\psi^0)$, we have 
$$
\align
&\Gam_{\neg_X}\neg(E^0_{\ss-style X })= \{\, (i_v \phi^0, i_v \psi^0 ) 
\,|\, v \in \Gam_{\neg_X}\neg( TX)\, \}\\
&\Gam_{\neg_X}\neg(E^1_{\ss-style X })= \{ \,(\theta\w i_v \phi^0,\theta\w i_v\psi^0
)\,|\,\theta\in
\Gam_{\neg_X}\neg(\w^1),  v \in \Gam_{\neg_X}\neg(TX)\,\} \\
&\Gam_{\neg_X}\neg(E^2_{\ss-style X }) =\{\, (\theta\w\a,\theta\w\b )\,|\, 
\theta \in \Gam_{\neg_X}\neg(\w^1),(\a,\b) \in
\Gam_{\neg_X}\neg(E^1_{\ss-style X })\, \}
\endalign
$$
Let $i_{_M}$ be the inclusion $M \to X$. 
Then by using the pull back $i^*_{_M}$ of differential forms in 
the second component of each $E^i_X$,
we obtain vector bundles $E_{_M}^i$ over $M$ $(i=0.1,2),$
$$
\align
&\Gam_{\neg_M}\neg(E^0_{\ss-style {M}})= \{\,i^*_{_M} (i_v \psi^0 ) 
\,|\, v \in \Gam_{\neg_X}\neg( TX)\, \}\\
&\Gam_{\neg_M}\neg(E^1_{\ss-style M })= \{ \,i_{_M}^*(\theta\w i_v\psi^0 )\,|
\,\theta\in \Gam_{\neg_M}\neg(\w^1), 
v \in \Gam_{\neg_X}\neg(TX)\,\} \\
&\Gam_{\neg_M}\neg(E^2_{\ss-style M } )=\{\, i_{_M}^*(\theta\w\b )\,|\, 
\theta \in \Gam_{\neg_M}\neg(\w^1),(\a,\b) \in \Gam_{\neg_X}\neg(E^1_{\ss-style X
})\, \}.
\endalign
$$
We denote by \#$_{_M}$ the complex on $M$, 
$$\CD
\Gam_{\neg_M}\neg(E^0_{\ss-style M })@>d_{M,0}>>\Gam_{\neg_M}\neg(E^1_{\ss-style M })@>d_{M,1}>>\Gam_{\neg_M}\neg(E^2_{\ss-style M })
\endCD
\tag\#$_{_M}$
$$
From our construction of $E^i_M$,
we have the map $\kappa$ in terms of the pull back $i_{_M}^*$,
$$
\kappa\: \Gam_{\neg_X}\neg(E^i_X) \arrow \Gam_{\neg_M}\neg(E^i_M).
$$
Then the map $\kappa$ is the map from the complex \#$_{_X}$ to 
the complex \#$_M$.
Hence we have the short exact sequence of complexes, 
$$
\CD 
0@.0@.0 \\
@VVV@VVV@VVV\\
\Gam_{\neg_{X,M}}\neg(E^0_{\ss-style X })
@>>>\Gam_{\neg_{X,M}}\neg(E^1_{\ss-style X })@>>>\Gam_{\neg_{X,M}}\neg(E^2_{\ss-style X })\\
@VVV@VVV@VVV\\
\Gam_{\neg_X}\neg(E^0_{\ss-style X })@>d_{X,0}>>\Gam_{\neg_X}\neg(E^1_{\ss-style X })@>d_{X,1}>>\Gam_{\neg_X}\neg(E^2_{\ss-style X })\\
@VVV@VVV@VV V\\
\Gam_{\neg_M}\neg ( E^0_{\ss-style M })@>d_{M,0}>>\Gam_{\neg_M}\neg( E^1_{\ss-style M })@>d_{M,1}>>\Gam_{\neg_M}\neg( E^2_{\ss-style M })\\
@VVV@VVV@VVV\\
0@.0@.0
\endCD
\tag 1$$
where $\Gam_{\neg\neg_{X,M}}\neg(E^*_{\ss-style X})$ is given as 
$$
\align
&\Gam_{\neg_{X,M}}\neg(E^0_{\ss-style X })= \{\, (i_v \phi^0, i_v \psi^0 ) \in \Gam_{\neg_X}\neg(E^0_{\ss-style X })
\,| \,i^*_{_M} (i_v\psi^0)=0 \,\}\\
&\Gam_{\neg_{X,M}}\neg(E^1_{\ss-style X })= \{ \,(\theta\w i_v \phi^0,\theta\w i_v\psi^0
)  \in \Gam_{\neg_X}\neg(E^1_{\ss-style X })\,|\,    i_{_M}^* (\theta\w i_v\psi^0)=0 \,\} \\
&\Gam_{\neg_{X,M}}\neg(E^2_{\ss-style X }) =\{\, (\theta\w\a,\theta\w\b )\in \Gam_{\neg_X}\neg(E^2_{\ss-style X })\,|\, 
i_{_M}^*(\theta\w\b)=0 \,
 \},
\endalign
$$ 
We denote by \#$_{X,M}$ 
the complex $(\Gam_{\neg_{X,M}}\neg(E^*_{_X}), d_{\ss-style{X,M}})$,
$$\CD
\Gam_{\neg_{X,M}}\neg(E^0_{\ss-style X })
@>d_{X,M}>>\Gam_{\neg_{X,M}}\neg(E^1_{\ss-style X
})@>d_{X,M}>>\Gam_{\neg_{X,M}}\neg(E^2_{\ss-style X })
\endCD\tag\#$_{\neg\scriptscriptstyle{X\neg, M}}$ 
$$ 
\proclaim{Theorem 4-1-5}
Let $\M_{_{\neg K\neg E}}\neg(X,M)$ be the relative moduli space. 
Then $\M_{_{\neg K\neg E}}\neg(X,M)$ is a smooth manifold 
(in particular, $\M_{_{\neg K\neg E}}\neg(X,M)$ is Hausdorff).
Further local coordinates of 
$\M_{_{\neg K\neg E}}\neg(X,M)$ is
given by 
an open ball of the cohomology group 
H$^1(\#_{\neg\scriptscriptstyle{X\neg,
M}})$, where H$^1(\#_{\neg\scriptscriptstyle{X\neg,M}})$ is the first cohomology
group of the complex 
\#$_{\neg\scriptscriptstyle{X\neg,M}}$. 
\endproclaim
We will prove theorem 4-1-5 in the rest of this section.
\proclaim{Lemma 4-1-6}
Let $H^i(\#_{_X})$ be the cohomology group of the complex \#$_{_X}$.
Then we have  
$$
H^0(\#_{_X}) = H^1(X).
$$
\endproclaim
\demo{Proof}
Let $i_v\Phi^0=(i_v\Ome^0, i_v\ome^0)$ be an element of 
$\Gam_{\neg_X}\neg(E^0_{ X})$,
for some real vector
$v \in \Gam_{\neg_X}\neg(TX)$. Let $g_{\Phi^0}$ be the
metric of
$X$  corresponding to $\Phi^0=(\Ome^0,\ome^0)$ and $*$ the Hodge star operator
with  respect to $g_{\Phi^0}$. Then we see 
that 
$$
\align 
&*i_v \Ome^0 =c_1\, i_v\ome^0 \w \ol{\Ome^0} \\
&*i_v\ome^0 = c_2 \,i_v\Ome^0 \w  \ol{\Ome^0}+\ol{c_2}\,\,i_v\ol{\Ome}^0 \w
\Ome^0,
\tag2\endalign
$$where $c_1, c_2$ are constants depending only on the dimension $n$.
Hence if $di_v \ome^0 =0$ and $dI_V\Ome^0 =$, the $i_v\ome^0$ is a harmonic
$1$ form with respect to the metric $g_{\Phi^0}$ from (1). 
We define the complex Hodge star operator $*_{_\C}$ by 
$$
*_{_\C}\: \w^{p,0}_{_X} \to \w^{n-p,0}_{_X},
$$
$$
c_2\,(*_{_\C} \a) \w\ol{\Ome}^0= * \a.
$$
Then we have 
$$
*_{_\C}i_v \ome^0 = i_v \Ome^0.
\tag3$$
Conversely any $1$ from is written as $i_v \ome^0$ since 
$\ome^0$ is a symplectic form. 
If $i_v\ome^0$ is a harmonic one form, then from (3) $i_v\Ome^0$ is also
harmonic.  Hence we have the result.  
\qed\enddemo
Then from lemma 4-1-6 we have 
\proclaim{Lemma 4-1-7}
$$H^0(\#_{_X})= \{ \, ( a,b)\in \H^{n-1}(X,\C)\oplus \H^1(X) \,|\, 
*_{_\C} a = (b)^{1,0}\,\},
$$ where $\H^i(X)$ denotes harmonic forms on $X$ with respect to $g_{\Phi^0}$.
\endproclaim
\demo{Proof}
This directly follows from lemma 6.
\enddemo

\proclaim{Lemma 4-1-8}
Let $N_{ \ss-style M}$ be the normal bundle on the special lagrangian
$M$ with respect to a K\"ahler-Einstein structure $\Phi^0=(\Ome^0,\ome^0)$. 
Then we have the canonical identification, 
$$
N_{\ss-style M} \cong E^0_{\ss-style M }
$$
\endproclaim
\demo{Proof}We consider a splitting
$\Gam_{\neg_M}\neg(N_{\ss-style M})\to
\Gam_{\neg_X}\neg(TX)$ of the exact sequence, 
$$
\CD
0@>>> \Gam_{\neg_{X,M}}\neg(TX)
@>>>\Gam_{\neg_X}\neg(TX)@>>>\Gam_{\neg_M}\neg(N_{\ss-style M})@>>>0.
\endCD
$$ 
Since $M$ is a special lagrangian, $i_{_M}\neg^*\psi^0=0$ for $\Phi^0
=(\phi^0,\psi^0)$,where $\psi=((\Ome^0)^{Im}, \ome^0)$.
Using the splitting we have the identification:
$$\align
&N_{\ss-style M} \cong E^0_{\ss-style M },\\
&\quad v \mapsto i_{_M}^* i_v \psi^0.
\endalign$$
The identification does not depend on a choice of a splitting. 
\qed\enddemo
\proclaim{Lemma 4-1-9}
$E^0_{\ss-style M }$ is the set of self dual forms of $\w_M^1\oplus \w_M^{n-1}$, 
$$
E^0_{\ss-style M } =\{\, ( a^1,a^{n-1}) \in \w_M^1\oplus \w_M^{n-1}\, |\, 
*_{_M}\neg a^1 = a^{n-1}\, \},
$$
where $*_{_M}$ is the Hodge star operator with respect to the pull
back metric 
$i_{_M}^*g_{\Phi^0}$ on $M$.
\endproclaim
\demo{Proof}
Any element of $E^0_{\ss-style M }$ is given as 
$i_{_M}^*(i_v\psi^0)$. Then we see that 
$$
i_{_M}^*(i_v\Ome)^{Im}= *_{_M}  i_{_M}^*(i_v\ome^0).
\tag4
$$
Hence $i_{_M}^*(i_v\psi^0)$ is a self-dual form.
Conversely any self dual from is written as $i_{_M}^*(i_v\psi^0)$ for 
some $v \in N_{\ss-style M}$.
\qed\enddemo
\proclaim{Lemma 4-1-10}
Let H$^0(\#_{_M})$ be the cohomology group of the complex \#$_M$. 
Then we have 
$$
H^0(\#_{_M}) = H^1(M).
$$
\endproclaim
\demo{Proof}
From lemma 4-1-9, H$^0(\#_{_M})$ is the set of self dual harmonic forms
$(\H^1(M)\oplus\H^{n-1}(X))^+$. Hence we have the result. 
\qed\enddemo
\proclaim{Lemma 4-1-11}The image
$d_{\ss-style{M,0}}\Gam_{\neg_M}\neg(E^0_{\ss-style M })$
coincides with  the closed subset $d \Gam_{\neg_M}\neg(\w^1_M\oplus\w^{n-1}_M)$.
\endproclaim
\demo{Proof}
This follows from lemma 4-1-9.
\qed\enddemo
Let $\a$ be a real harmonic one form on $X$. Then we see that $*_{_\C}
(\a)^{1,0}
\in 
\Gam_{\neg_X}\neg( \w^{n-1,0})$ 
is also Harmonic. Hence by using the pull back $i_{_M}^*$, from (3) and (4) 
we have  
$$
i_{_M}^* (*_{_\C} \a^{1,0} )^{Im}=*_{_M}i_{_M}^*\a.
$$ 
Hence we have the commutative diagram:
\proclaim{Lemma 4-1-12}
$$
\CD
H^1(X) @>\gam_{H^1}>>H^1(M)\\
@V{*_{_\C}}VV @V{*_{_M}} VV \\
H^{n-1}(X,\C)@>>>H^{n-1}(M,\R),
\endCD
$$
where $\gam_{H^1}$ is the induced map from the pull back $i^*_{_M}$.
\endproclaim 
\proclaim{Lemma 4-1-13}
Let $H^1(\#_{_M})$ be the first cohomology group of the complex \#$_M$. 
Then we have 
$$
H^1(\#_{_M}) \cong H^n(M,\R)\oplus H^2(M,\R).
$$
\endproclaim
\demo{Proof}Let $I=I_{\Ome^0}$ be the complex structure on $X$
corresponding to 
$\Ome^0$. Since $I$ is an element of End $(TX)$, we have 
$\rho_{I}(\Ome^0) \in \Gam_{\neg_X}\neg(E^1_{\ss-style X })$ as in section
one.
We see that 
$\rho_{I}\Ome^0 =i\Ome^0$ since $\Ome^0 \in \w^{n,0}$. 
Since $M$ is a special lagrangian, $i_{_M}^*(\Ome^0)^{Re} = $vol$_{_M}$.
Hence 
$$i_{_M}^* ( (\rho_{I}\Ome^0)^{Im}) =i_{_M}^*(-\Ome^{Re})= -vol_M.
\tag5$$
Since $\ome$ is a K\"ahler form, $\rho_{I}\ome^0 = \ome^0$.
Hence $( vol_M,i_{_M}^*\ome)=(vol_M,0)$ is an element of $\Gam_{\neg_M}\neg(E^0_{\ss-style M })$. 
We also see that $\{\,i_v\ome\,|\,v \in \Gam_{\neg_X}\neg(TX)\,\} 
=\Gam_{\neg_X}\neg(\w^2_X)$. Hence 
$\{\,i_{_M}^*(i_v\ome)\,|\,v \in
\Gam_{\neg_X}\neg(TX)\,\}=\Gam_{\neg_M}\neg( \w^2_M)$. Then from (5) we
see that $\Gam_{\neg_M}\neg( E^1_{\ss-style M })= \Gam_{\neg_M}\neg(\w^n_M\oplus\w^2_M)$. Then from lemma 11 we
have the result.
\qed\enddemo
Let $\gam_{\neg_{H^2}}\:H^2(X,\R)\to H^2(M,\R)$ be the induced map 
from the pull back $i^*_{_M}$. We denote by Image $\gam_{\neg_{H^2}}$ the image of
the map $\gam_{\neg_{H^2}}$.
As in (1) we have the map from the complex \#$_{_X}$ to
the complex
$\#_{_M}$. Then we have the map of cohomologies,
$$
\gam^1_{_{\#}} \: H^1(\#_{_X}) \to H^1(\#_{_M}).
$$
\proclaim{Lemma 4-1-14}
The image of the map $H^1(\#_{_X}) \to H^1(\#_{_M})$ is the direct sum 
H$^n(M,\R)\oplus $Image $\gam_{\ss-style{H^2}}$, so that is, we have
the commutative diagram: 
$$\xy
 (-20,0)*+{H^1(\#_{\ss-style X})}="a",
 (20,0)*+{H^n(M,\R) \oplus \text{Image }\gam_{\ss-style H^2} }="b",
 (20,-20)*+{H^1(\#_{\ss-style M})}="c",
(60,0)*+{0}="d"
\ar "a";"b" 
\ar "b";"c"
\ar_{\gam^1_{\ss-style {\#}}} "a";"c"
\ar "b";"d"
  \endxy$$
\endproclaim
\demo{Proof}Since $\rho_{I}(\Ome^0) = i\Ome^0 $ and
$\rho_I{\ome^0}=\ome$, $(i\Ome^0,\ome)$ is an element of
$\Gam_{\neg_X}\neg(E^1_{\ss-style X })$. Hence $i_{_M}^* ( (i\Ome^0)^{Im},\ome)= (-vol_M,0)$ 
is an element of $\Gam_{\neg_M}\neg(E^1_{\ss-style M })$. From proposition 3-2-8 in section 3, 
$P_{\ss-style{ H^2}}\: H^1(\#_{_X}) \to H^2(X,\R)$ is surjective. Hence
we have the result. 
\qed\enddemo
We also denote by $\gam_{\neg_{H^1}}$ the induced map from $i_{_M}^*$,
$\gam_{\ss-style{H^1}}\:H^1(X,\R)\to H^1(M,\R)$.
From (1)
we have the long exact sequence:
$$\minCDarrowwidth{0.3cm}
\CD
H^0(\#_{_X})@>>>H^0(\#_{_M})@>>>H^1(
\#_{\neg\ss-style{X\neg,M}})@>>>H^1(\#_{_X})@>\gam^1_{_\#}>>\neg\neg\neg\neg\neg\neg\neg\neg
H^1(\#_{_M})\\ @| @|@.@.@|\\
H^1(X)@>\gam_{\neg_{H^1}}>>H^1(M)@.@.@.H^n(M)\oplus H^2(M)
\endCD$$
Hence we have an exact sequence:
$$\minCDarrowwidth{0.5cm}
\CD
0@>>>\text{Coker
}\gam_{_{H^1}}@>>>H^1(\#_{\neg\ss-style{X\neg,M}})@>>>\text {Ker
}\gam^1_{_\#}@>>>0.
\endCD\tag6$$
\proclaim{Proposition 4-1-15}
The dimension of H$^1(\#_{\neg\ss-style{X\neg,M}})$ is invariant under
local deformations of 
$\Phi^0 \in \til{\M}_{_{\neg K\neg E}}\neg(X,M)$.
\endproclaim
\demo{Proof}
From the exact sequence (6), 
$$\dim H^1(\#_{\neg\ss-style{X\neg,M}}) =\dim\text{ Coker }
\gam_{_{H^1}} + \dim \text{ Ker } \gam^1_{_\#}.
$$
From lemma 14 ker $\gam^1_{_\#}$ is the ker of the map H$^1(\#_{_X}) \to
H^n(M)\oplus 
$Image $\gam_{\neg_{H^2}}$, where Image $\gam_{\neg_{H^2}}$ does not depend on
$\Phi^0$. Hence $\dim$ Ker
$\gam^1_{_\#}$ is invariant since H$^1(\#_{_X})$  is invariant from proposition
1-5 in section one. Hence we see that
H$^1(\#_{\neg\ss-style{X\neg,M}})$ is also invariant  since image
$\gam_{\neg_{H^1}}$ is topological.
\qed\enddemo
As in section 3, the complex \#$_{_X}$ is a subcomplex of the following de Rham
complex dR$_{\ss-style X}$: 
$$\minCDarrowwidth{0.2cm}\CD
\cdots@>>>\Gam_{\neg_X}\neg (
\w^{\ss-style{n-1}}_{\ss-style\C}\oplus\w^{\ss-style 1}
)@>d_X>>\Gam_{\neg_X}\neg(\w^{\ss-style n}_{\ss-style \C}\oplus
\w{\ss-style^2})@>d_X>>
\Gam_{\neg_X}\neg(\w^{\ss-style{n+1}}_{\ss-style\C}\oplus
\w^{ \ss-style 3})@>>>\cdots.
\endCD
$$
We also see that the complex $\#_{_M}$ is a subcomplex of de Rham complex 
dR$_{\ss-style M}$:
$$
\minCDarrowwidth{0.3cm}
\CD
\cdots@>>>\Gam_{\neg_M}\neg(\w^{n-1}\oplus\w^1)@>d_M>>\Gam_{\neg_M}\neg(\w^n\oplus\w^2)@>d_M>>\Gam_{\neg_M}\neg(\w^{n+1}
\oplus \w^3 )@>>>\cdots.
\endCD
$$
Hence we see that the complex $\#_{\neg\ss-style{X\neg,M}}$ is a
subcomplex of the following  relative de Rham complex
dR$_{\ss-style{X\neg,M}}$: 
$$\minCDarrowwidth{0.2cm}
\CD
\cdots@>>>\Gam_{\neg_{X,M}}\neg(\w^{n-1}_\C\oplus\w^1 )@>d>>\Gam_{\neg_{X,M}}\neg(\w^n_\C\oplus
\w^2)@>d>>
\Gam_{\neg_{X,M}}\neg(\w^{n+1}_\C\oplus \w^3)@>>>\cdots,
\endCD
$$
where $\Gam_{\neg_{X,M}}\neg(\w^{p}_\C\oplus\w^q)$ is given as 
$$
\Gam_{\neg_{X,M}}\neg(\w^p_\C\oplus\w^q) =\{ \, ( a,b) \in \Gam_{\neg_X}\neg(\w^p_\C\oplus\w^q) 
\,|\, i^*_{_M}a^{Im}=0, \, i^*_{_M} b=0 \,\}.
$$
Then we have the map between short exact sequences of complexes: 
$$\CD
0@>>>\#_{\neg_{X,M}}@>>>\#_{_X}@>>>\#_{_M}@>>>0 \\
@.@VVV @VVV @VVV@.\\
0@>>>\text{dR}_{\neg_{X,M}}@>>>\text{dR}_{_X}@>>>\text{dR}_{_M}@>>>0
\endCD$$
Hence we have the map between long exact sequences: 
$$\minCDarrowwidth{0.5cm}\CD
H^0(\#_{_X})@>\gam^0_{_{\#}}>>H^0(\#_{_M})@>>>H^1(\#_{\neg\ss-style{X\neg,M}})@>>>H^1(\#_{_X})\\
@VVV @VVV@VVp_{_{\neg\ss-style{X\neg,M}}}V@VVp_{_X}V\\
H^0(dR_{_X})@>\gam_{_{dR}}>>H^0(dR_{_M})@>>>H^1(dR_{\neg\ss-style{X\neg,M}})@>>>
H^1(dR_{_X})
\endCD
$$
\proclaim{Proposition 4-1-16}
The map $p_{\neg\ss-style{X\neg,M}}\:$
H$^1(\#_{\neg\ss-style{X\neg,M}}) \to
$H$^1($dR$_{\neg\ss-style{X\neg,M}})$ is injective.
\endproclaim
\demo{Proof}
We consider the induced map
between quotient vector
spaces in the diagram : 
$$
H^0(\#_{_M}) /\gam^0_{_\#}H^0(\#_{_X}) \to H^0(dR_{_M}) /\gam_{_{dR}}H^0(
dR_{_X}).
\tag7
$$
From lemma 4-1-12 we see the map (7) is injective. From proposition 3-2-3 of section 3, 
$p_X\: H^1(\#_{\ss-style X})\to H^1( dR_{\ss-style {X}})$ is also injective. Hence we
have the result. 
\qed\enddemo
\subhead \S 4-2
\endsubhead
We shall construct a slice on
$\til{\M}(X,M)$ for the action of 
$\Diff_0(X)$. At first we consider the slice $S_{\Phi^0}=S_{\Phi^0}(X)$ 
in section
$2$ and define the map $\lam$: 
$$\align
\lam\: &S_{\Phi^0}(X) \arrow H^n(M,\R)\oplus H^2(M,\R),\\
&\Phi=(\phi,\psi) \mapsto [ i_{_{\ss-style  M}}^*\psi]_{dR},
\endalign$$
where $\Phi\in S_{\Phi^0}$.
\proclaim{Lemma 4-2-1}
We define $\h{S}_{\Phi^0}(X)$ by 
$$
\h{S}_{\Phi^0}(X) =\{\, \Phi \in S_{\Phi^0}(X) \, |\, 
\lam(\Phi)=0 \, \}.
$$
Then a neighborhood of $\h{S}_{\Phi^0}(X)$ at $\Phi^0$ is homeomorphic to 
an open set of Ker $\gam^1_{\ss-style\#}$ in $H^1(\#_{{\ss-style X}})$, So that is,
$\h{S}_{\Phi^0}(X)$ is a manifold around $\Phi^0$ with,
$$
T_{\Phi^0}\h{S}_{\Phi^0}(X) = \text{ Ker }\gam^1_{_\#}.
$$
\endproclaim
\demo{Proof}
From lemma 4-1-13 the image of the map $\gam^1_{\#}\:H^1(\#_{\ss-style X}) \to
H^n(M)\oplus H^2(M)$ is  given by $H^n(M)\oplus\text{Image }\gam_{_{H^2}}$. Hence
Image $\gam^1_{_\#}$ does not depend on 
$\Phi^0 \in \wtil{\M}_{KE}(X,M)$. Let $\Phi $ be an element of 
the slice $S_{\Phi^0}(X)$. Then we have the map from the tangent space T$_\Phi
S_{\Phi^0}(X) $ to $H^n(M)\oplus H^2(M)$ by taking cohomology classes.
Since $\Phi$ is a topological calibration, the image of this map 
coincides $H^n(M)\oplus\text{Image }\gam_{_{H^2}}$. 
(Note that the difference between T$_{\Phi}S_{\Phi}$ and T$_{\Phi}S_{\Phi^0}$
is given by exact forms.)
Hence we see that 
the image $\lam( S_{\Phi^0})$ is $H^n(M)\oplus\text{Image }\gam_{_{H^2}}$.
By using implicit function theorem we see that the inverse image $\lam^{-1}(0)$
is a manifold around $\Phi^0$
\qed\enddemo
As in section 2, we define the Sobolev space $L^2_s( M,\w^p)$, 
for $s > k + \frac {2n}2$, where $s$ is sufficiently large.
 Then we denote by 
$L^2_s(E^i_{\ss-style M})=L^2_s( M, E^i_{\ss-style M})$ 
the Sobolev space $L^2_s(M,\w^{n-1+i}\oplus\w^{i+1})\cap 
C^k (M,E^i_{\ss-style M})$ for $i=0,1,2$. 
Since the complex \#$_{\ss-style  M}$ is  an elliptic complex, we have the Hodge
decomposition, 
$$
L^2_{s}(M,E^0_{\ss-style M}) = \H^0(\#_{\ss-style  M}) \oplus d^*_{_{M,0}}L^2_{s+1}(M,E^1_{\ss-style M}),
$$
where $\H^1(\#_{\ss-style  M})$ is the Ker $d_{_{M,0}}$. 
We have the commutative diagram: 
$$
\CD
L^2_{s}(X, E^0_{\ss-style X}) @>d_{X,0}>>L^2_{s-1}(X, E^1_{\ss-style X})@>>>\\
@V\kappa^0  VV @V\kappa^1 VV\\
L^2_s(M,E^0_{\ss-style M}) @>d_{M,0}>>L^2_{s-1}(X, E^1_{\ss-style M})@>>>.
\endCD
$$
Then we denote by $\gam^0=\gam^0_{\ss-style \#}$ the induced map on the cohomologies
from
$\kappa^0$, 
$$
\gam^0\: H^0(\#_{\ss-style X}) \arrow H^0(\#_{\ss-style  M}).
$$

Let $\gam^0(\text{Ker }d_{_{X,0}})=\gam^0(\H^0(\#_{\ss-style X}))$ be the image of $\H^0(\#_{\ss-style X})$ 
in $L^2_s(E^0_{\ss-style M})$. We denote by $(\text{Image }\gam^0)^\perp$ the 
orthogonal complement of the image $\gam^0(\H^0(\#_{\ss-style X}))$ in 
$L^2_s(M,E^0_{\ss-style M})$. Since $\text{ Image }\gam^0$ is the subspace of 
$\H^0(\#_{\ss-style  M})$, we have the decomposition: 
$$\align
L^2_s( M, E^0_{\ss-style M})  
= &\text{Image }\gam^0\oplus (\text{Image }\gam^0)^\perp\\
=& \text{Image }\gam^0\oplus \H^0 (\#_{\ss-style X}, \#_{\ss-style  M} ) \oplus 
d^*_{_{M,0}}L^2_{s+1}(M, E^1_{\ss-style M}), 
\endalign
$$
where $H^0(\#_{\ss-style X},\#_{\ss-style  M})$ denotes the orthogonal complement of 
Image $\gam^0$ in $ H^0(\#_{\ss-style  M})$. We recall the identification $\Gam(E^0_{\ss-style M}) \cong
\Gam(N_{\ss-style M})$ in lemma 8 and  fix a splitting $\Gam_{\ss-style M}( N_{\ss-style M}) \to \Gam_{\ss-style X}(TX)$ of the
exact sequence,
$$\CD
0@>>>\Gam_{\ss-style{X,M}}(TX)@>>>\Gam_{\ss-style
X}(TX)@>>>\Gam_{\ss-style M}(N_{\ss-style M})@>>>0.
\endCD
$$
Under the identification $\Gam(E^0_{\ss-style M}) \cong \Gam(N_{\ss-style M})$, 
$v\in L^2_s( E^0_{\ss-style M})\cong L^2_s(N_{\ss-style M})$ is regarded with an element of $L^2_s(TX)$ by 
the splitting map. Let $\exp$ be the exponential map with respect to the metric
$g_{\Phi^0}$.
We define the map $F_{_{\ss-style {X,M}}}$ by 
$$\align
F_{_{\ss-style {X,M}}} \: (&\text{Image }\gam^0)^\perp\times \h{S}_{\Phi^0} \arrow dL^2_s (
M,E^0_{\ss-style M}),\\
&\qquad( \quad v\, \, ,\,\quad\Phi \quad) \arrow \quad i_{\ss-style M}^* \exp_v^*
\psi,
\endalign
$$
where $\Phi =(\phi,\psi)\in\h{S}_{\Phi^0}$ and 
$\exp_v^*\psi=(\,\exp_v^* (\Ome)^{Im}, \exp_v^*\ome \,)$ is the pull back 
of $\psi$ by $\exp_v \in \Diff_0(X)$. 
From Lemma 4-2-1 $i_{_M}^*\exp_v^*\psi$ is an
exact form and we see that the image of the map $F_{\ss-style{X,M}}$
is $dL^2_s(M,E^0_{\ss-style M})$ from lemma 4-1-1.
\proclaim{Proposition 4-2-2} 
The inverse image $F^{-1}_{_{\ss-style {X,M}}}(0)$ is a manifold around $(0,\Phi^0)$. 
\endproclaim
\demo{Proof}
The differential of $F_{_{\ss-style {X,M}}}$ at $(0, \Phi^0)$ is given by 
$$
dF_{_{\ss-style {X,M}}}( \dot{v},\dot{\Phi}) \arrow 
i_{_{\ss-style  M}}^* d\,i_{\dot{v}}\psi^0 
+ i_{_{\ss-style  M}}^* \dot{\psi},
$$where $(\dot{v},\dot{\Phi})\in (\text{Image
}\gam^0)^\perp\oplus  T_{\Phi^0 }\h{S}_{\Phi^0}$ 
and $\dot{\Phi}=(\dot{\phi},\dot{\psi})$. Hence we see that
$dF_{_{\ss-style {X,M}}}$ is surjective.  From the implicit function theorem, we
have the result.
\qed\enddemo
\proclaim{Definition 4-2-3} 
Let $\exp$ be the exponential map corresponding to the metric $g_{\Phi^0}$.
We define the slice $S_{\Phi^0}(X,M)$ by 
$$
S_{\Phi^0}(X,M) = \{ \, \exp_v^*\Phi \in \wtil{\M}(X)\cap  U\, |
\, (v, \Phi ) \in F_{\ss-style {X,M}}^{-1}(0)\, \},
$$
where $U$ is a sufficiently small neighborhood of $\wtil{\M}(X)$ at $\Phi^0$.
\endproclaim
Then we have 
\proclaim{Proposition 4-2-4}
The slice $S_{\Phi^0}(X,M)$ is a submanifold of $\wtil{\M}(X,M)$.
\endproclaim
\demo{Proof}
From elliptic regularity we see that $(v,\Phi) \in F^{-1}(0)_{_{\ss-style {X,M}}}$ is 
smooth. 
From definition, $i_{\ss-style M}^* \exp_v^*\Phi =0$. Hence $S_{\Phi^0}(X,M) 
$ is a subset of $\wtil{\M}(X,M)$.  As in section 2, 
$U \cong V \times S_{\Phi^0}$, where $V$ is a neighborhood of 
\Diff$_0(X)$ at the identity. Hence $S_{\Phi^0}(X,M)$ is homeomorphic to 
a neighborhood of $F^{-1}_{_{\ss-style {X,M}}}(0)$. From Proposition 4-2-2 we have
the result.
\qed\enddemo
We recall the exact sequence (1):
$$
\CD
0@>>>\#_{\ss-style{X,M}}@>>>\#_{\ss-style X}@>\kappa>>
\#_{\ss-style M}@>>>0.
\endCD
$$
From the map $\kappa \:\#_{\ss-style  X} \to \#_{\ss-style  M}$, we define the following relative
complex:
$$
C^{\ss-style p}( \#_{\ss-style  X} \to \#_{\ss-style  M}) = \Gam_{\ss-style X}(
E^{\ss-style p}_{\ss-style X})
\oplus\Gam_{\ss-style M}(E^{\ss-style {p-1}}_{\ss-style  M}), 
$$
where the coboundary map $\del_{_{\ss-style {X,M}}}$ is  
$$
\del_{_{\ss-style {X,M}}}( a_{\ss-style X}^{\ss-style p}, b_{\ss-style M}^{\ss-style  p-1} ) = 
( d a_{\ss-style X}^{\ss-style p}, (-1)^{\ss-style  {p-1}} \kappa(a) + 
d b_{\ss-style M}^{\ss-style {p-1}} ),
$$where $(a^{\ss-style p}_{\ss-style X},b_{\ss-style M}^{\ss-style {p-1}})\in
\Gam_{\ss-style X}(E^{\ss-style p}_{\ss-style X})\oplus
\Gam_{\ss-style M}(E^{\ss-style {p-1}}_{\ss-style  M})$. ( Note that
$C^0(\#_{\ss-style  X}\to \#_{\ss-style  M})= \Gam_{\ss-style X}(E^0)$. ) Then we
have a cohomology group
$H^{\ss-style p}( \#_{\ss-style  X} \to \#_{\ss-style  M})$ of the complex $(
C^{\ss-style p} (
\#_{\ss-style  X}
\to
\#_{\ss-style  M} ) , 
\del_{\ss-style {X,M}}) $. 
\proclaim {Lemma 4-2-5}Let $H^p(\#_{_{\ss-style {X,M}}})$ be the cohomology group of
the complex 
$\#_{_{\ss-style {X,M}}}$. Then we have 
$$H^p(\#_{_{\ss-style {X,M}}}) \cong H^p(\#_{\ss-style  X}\to \#_{\ss-style  M}).$$
\endproclaim
\demo{Proof}Let $a_{\ss-style X}^{\ss-style  p}$ be an element of
$\Gam_{\ss-style {X,M}}(E^{ p}_{\ss-style X})$.  Then we define the map 
from $\Gam_{\ss-style {X,M}}(E^{ p}_{\ss-style X})$ to 
$C^{\ss-style p} (
\#_{\ss-style  X}
\to
\#_{\ss-style  M} ) $ by 
$$a_{\ss-style X}^{\ss-style  p}\mapsto ( a_{\ss-style X}^{\ss-style  p},0).$$  
We see that this map is quasi isomorphism between complexes.
\qed\enddemo
\proclaim{Proposition 4-2-6}
Let $S_{\Phi^0}(X,M)$ be the slice as in definition 4-2-3. Then the tangent space
of the slice is given by the relative cohomology group:
$$
T_{\Phi^0}S_{\Phi^0}(X,M) = H^1(\#_{\ss-style {X,M}}).
$$
\endproclaim
\demo{Proof}
From Proposition 4-2-4 the tangent space $T_{\Phi^0}S_{\Phi^0}(X,M)$ is ker
$dF_{\ss-style {X,M}}$. From Proposition 4-2-2 Ker $dF_{\ss-style {X,M}}$ is given as 
$$
\text{ Ker }  dF_{\ss-style {X,M}} = \left\{ \, ( \dot v, \dot \Phi ) \in
(\text{Image }\gam^0)^\perp 
\oplus T_{\Phi^0}\h{S}_{\Phi^0}\,\, \big|\, \,
i_{_{\ss-style  M}}^* d\,i_{\dot v} \psi^0 + i_{_{\ss-style  M}}^* \dot \psi=0\,
\right\}.
$$
Then $(\dot \Phi, i_{\dot v}\psi^0)\in $  Ker $dF_{\ss-style {X,M}}$ is 
a representative of $H^1(\#_{\ss-style  X}\to \#_{\ss-style  M})$. 
Hence we have the map Ker $dF_{\ss-style {X,M}}\to H^1(\#_{\ss-style  X}\to \#_{\ss-style  M})$.
From Hodge decomposition, we see that this map is 
bijective. 
\qed\enddemo
\proclaim{Proposition 4-2-7}
Let $\pi_{{\ss-style {X,M}}}\: \wtil{\M}(X,M)\to \M(X,M)$ be the natural projection.
We restrict $\pi_{{\ss-style {X,M}}}$ to the slice $S_{\Phi^0}(X,M)$. Then the
restricted map is injective and its image is an open set of $\M(X,M)$. 
\endproclaim
\demo{Proof}
Let $\Phi$ be an element of the slice $S_{\Phi^0}(X,M)$. Then $\Phi$ is a closed 
form of $\Gam_{\ss-style {X,M}}(E^1_{\ss-style X})$. 
Let $P$ be the map from the slice $S_{\Phi^0}(X,M)$ to the relative de
Rham 
cohomology group H$^1( dR_{\ss-style {X,M}})$, 
$$
P\:S_{\Phi^0}(X,M)\arrow H^1(dR_{\ss-style {X,M}}),
$$
$$
\Phi \arrow [\Phi]_{dR_{\ss-style{X,M}}}.
$$
From proposition 4-2-6, the differential of the map $P$ at $\Phi^0$ 
is an isomorphism 
$$
T_{\Phi^0}S_{\Phi^0}(X,M) \cong H^1(\#_{\ss-style {X,M}}).
$$
From proposition 4-1-16 $H^1(\#_{\ss-style {X,M}}) \to H^1(dR_{\ss-style {X,M}})$ is
injective.  Hence the map $P$ is injective if we take $S_{\Phi^0}(X,M)$ sufficiently
small. A class of $H^1( dR_{\ss-style {X,M}})$ is invariant under the action 
of Diff$_0(X,M)$. Hence the restricted map
$\pi_{\text{\fiverm X,M}}|_{S_{\Phi^0}(X,M)}$ is injective. let $U$ be a small
open set of $\wtil{\M}(X)$ at $\Phi^0$.  Then as in section 2, $U\cong V \times
S_{\Phi^0}$, where 
$V$ is a neighborhood of Diff$_0(X)$ at the identity. 
Let $\Phi $ be an element of $\wtil{\M}(X,M)\cap U$. 
Then $\Phi$ is written as $\Phi = f^*\h{\Phi}$, 
where $f \in V, \h{\Phi}\in S_{\Phi^0}$. 
Since $i_{\ss-style M}^*\Phi=0$, we have $[i_{\ss-style M}^* \h{\Phi}]_{dR} =0$. 
Hence $\h{\Phi}\in \h{S}_{\Phi^0}$. 
The open set $V$ of Diff$_0(X)$ at the identity is 
identified with an open set of $\Gam_{\ss-style X}(TX)$ at zero section. 
By using a splitting, we have the decomposition 
$$
\Gam_{\ss-style X}(TX) = \Gam_{\ss-style {X,M}}(TX)\oplus \Gam(N_{\ss-style M}).
$$
Let $u$ be an element of $\gam^0(H^0(\#_{\ss-style  X}))$. Then we see $ d\,i_u \Phi^0=0$. 
Since $U$ is sufficiently small, $\Phi$ is written as 
$$
\Phi = f_1^* \exp_v^* \h{\Phi},
$$
where $f_1 \in $Diff$_0(X,M)$ and $v \in (\text{Image }\gam)^\perp 
\subset H^0(\#_{\ss-style  M})$.
Hence Image $\pi_{\ss-style {X,M}}( S_{\Phi^0}(X,M) ) = \pi (U)$ is an open set.
\enddemo
We define the map $P_{\ss-style {X,M}}$ by 
$$
P_{\ss-style {X,M}} \: \M(X,M) \arrow H^1(dR_{\ss-style {X,M}}),
$$
$$
\Phi \arrow [\Phi]_{\ss-style {dR}}
$$
\proclaim{Theorem 4-2-8} 
The map $P$ is locally injective.
\endproclaim
\demo{Proof}
This follows from proposition 4-1-16.
\qed\enddemo
\demo{Proof of theorem 4-1-5} 
Let $\wtil{\M}_{\ss-style {KE}}(X,M)$ be as in definition 1. 
Since $\wtil{\M}_{\ss-style {KE}}(X,M)$ is a submanifold of 
$\wtil{\M}_{\ss-style{KE}}(X)$, we have the distance $d$ on 
$\wtil{\M}_{\ss-style {KE}}(X,M)$ as in proposition 2-9. 
Then from theorem 4-2-8 we see that 
$d$ gives a distance on the moduli space $\M_{\ss-style{KE}}(X,M)$ 
as in a proof of
proposition 2-9. Hence $\M_{\ss-style{KE}}(X,M)$ is Hausdorff.
Let $f$ be an element of Diff$_0(X,M)$. Then $f$ gives 
a diffeomorphism from the slice 
$S_{\ss-style {\Phi^0}}(X,M) $ to $S_{\ss-style {f^*{\Phi^0}}}(X,M)$. 
Hence from proposition 4-2-7, we see that each slice is local coordinates of 
$\M(X,M)$. From Proposition 4-2-4 and 4-2-6, we have the result.
\qed\enddemo
We define $\widehat{\M}(X,M)$ by 
$$
\widehat{\M}(X,M)=\left\{\, \Phi\in \E^1_{KE}(X)\, |\, d\Phi =0,\,\, 
[i_{\ss-style M}^* \psi ]_{dR} =0\, \right\}/ \Diff_0(X).
$$
Then we have the natural map 
$$\align
{\M}&(X,M) \arrow \widehat{\M}(X,M),\\
&[\Phi]_{{\text{Diff}_0(X,M)}} \mapsto [\Phi]_{\text{Diff}_0(X)},
\endalign$$
where $[\Phi]_{\text{Diff}_0(X,M)}$ ( resp. $  [\Phi]_{\text{Diff}_0(X)}$)
denotes a class of $\M(X,M)$ (resp. $\wtil{\M}(X,M)$ ).
We denote by $\widehat{\M}_0(X,M)$ the image of this map.
\proclaim{proposition 4-2-9}
 $\widehat{\M}_0(X,M)$ is a smooth manifold and local coordinates are
given  by an open ball of Ker $\gam^1_{\ss-style\#} \subset H^1(\#_{\ss-style  X})$
\endproclaim
\demo{Proof}
Local coordinates of $\widehat{\M}_0(X,M)$ is given by $\h{S}_{\Phi^0}$ 
as in Lemma 4-2-1. 
\enddemo
\proclaim{Theorem 4-2-10}
${\M}(X,M) \arrow \widehat{\M}_0(X,M)$ is a fibre bundle and 
local coordinates of each fibre is given by an open ball of 
 Coker $\gam_{\ss-style{H^1}}$.
\endproclaim
\demo{Proof}
Let $\Phi^0$ be an element of $\wtil{\M}(X,M)$. Then the fibre through 
$\Phi^0$ is written as the quotient,
$$
\{ \, f^*\Phi^0 \,|\, i_{\ss-style M}^*f^*\psi^0 =0, f \in \Diff_0(X) \,\}
/\{\, g^*\Phi^0 \, |\, g \in \Diff_0(X,M)\, \}.
$$
An neighborhood of Diff$_0(X)/\Diff_0(X,M)$ is given by 
an open set of $\Gam_{\ss-style M}(N_{\ss-style M})$. Under the identification $\Gam_{\ss-style M}(N_{\ss-style M}) \cong 
\Gam_{\ss-style M}(E^0_{\ss-style M})$, we see that each fibre is parameterized by 
an open ball of coker $\gam_{\ss-style {H^1}} \subset H^0(\#_{\ss-style  M})$. 
\qed\enddemo
Each fibre can be regarded as the moduli space of special lagrangian 
submanifolds with respect to a fixed K\"ahler-Einstein structure $\Phi$. 
\proclaim{Theorem 4-2-11}
Let $\wtil{\M}_{\ss-style{KE}}(X,M)$ be as in definition 4-1-1. 
We denote by $\Diff(X,M)$ the group of diffeomprphisims of $X$ preserving $M$.
There is the action of $\Diff(X,M)$ on $\wtil{\M}_{\ss-style{KE}}(X,M)$.
Then the quotient $\wtil{\M}_{\ss-style{KE}}(X,M)/\Diff(X,M)$ is an orbifold.
\endproclaim
\demo{Proof} 
The slice $S_{\Phi^0}(X,M)$ is local coordinates of
$\M_{\ss-style{KE}}(X,M)$ and  invariant under the action of $\Diff(X,M)$. 
Hence we see that the moduli space $\wtil{\M}_{\ss-style{KE}}(X,M)/\Diff(X,M)$ is locally
homoemorphic to  the quotient space $S_{\Phi^0}(X,M)/I_{\ss-style{X.M}}$, where
$I_{\ss-style{X.M}}$ is the isotropy, 
$$
I_{\ss-style{X.M}} = \{ \, f \in \Diff(X,M)\,| f^* \Phi^0 =\Phi^0 \, \}.
$$
As in proof of proposition 2-9, $\Diff(X,M)$ acts on $S_{\Phi^0}(X,M)$ isometrically. 
Hence we see that there is an open set $V$ of $T_{\ss-style{\Phi^0}} S_{\ss-style{\Phi^0}}(X,M)$
with the action of
$I_{\ss-style{X,M}}$  such that the quotient $V/I_{\ss-style{X.M}}$ is homeomorphic to
$S_{\Phi^0}(X,M)/I_{\ss-style{X.M}}$.  Since $T_{\Phi^0}S_{\Phi^0}(X,M)$ is isomorphic to
$H^1(\#_{\ss-style{X.M}})$ and the action of
$I_{\ss-style{X.M}}$ on $H^1(\#_{\ss-style{X.M}})$ is a isometry with respect to the
metric $g_{\ss-style {\Phi^0}}$. The action of
$I_{\ss-style{X.M}}$ preserves integral  cohomology class. Hence from proposition 4-1-16 we see
that 
$V/I_{\ss-style{X.M}}$ is  the quotient by a finite group. 
\qed\enddemo


\Refs
\widestnumber\key{10} 
\ref
\key 1
\by  A.L.Besse
\book Einstein manifolds
\publ  Ergebnisse der Mathematik und ihrer Grenzgebiete {\bf 10},
Springer-Verlag, Berlin-New York
\yr 1987
\endref
\ref
\key 2
\by P.Candelas and X.C.~de la Ossa
\paper  Moduli space of Calabi-Yau manifolds
\jour Nuclear Phys.  B 
\vol 355
\yr 1991
\pages 455--481
\endref
\ref 
\key 3
\by D.G.~Ebin
\paper The moduli space of riemannian metrics 
\jour Global Analysis, Proc. Symp. Pure Math. AMS 
\vol 15 
\yr 1968 
\pages 11-40
\endref
\ref 
\key 4
\by A.~Fujiki and G.~Schumacher 
\paper The moduli space of Extremal compact K\"ahler manifolds and
Generalized Weil-Perterson Metrics
\jour Publ.~RIMS,~Kyoto Univ
\vol 26.~No.1
\yr 1990 
\pages 101-183
\endref
\ref 
\key 5
\by R.~Goto
\paper 
Moduli spaces of topological
calibarations, Calabi-Yau, K\"ahler-Einstein, $\G$ 
and Spin$(7)$ structures
\paperinfo preprint 
\yr 2001 
\endref 
\ref 
\key 6 \by R.~Goto 
\paper On Hyper-K\"ahler manifolds of type $A_\infty$ and $D_\infty$
\jour Comm. Math. Phys
\vol 198
\yr 1998
\pages 469-491
\endref
\ref 
\key 7
\by F.~ R. ~Harvey
\book  Spinors and Calibrations
\publ  Academic Press, New York
\yr 1990
\endref
\ref \key 8
\by F.~ R.~ Harvey and H.~ B.~ Lawson
\paper Calibrated Geometries
\jour Acta Math
\vol 148 
\yr 1982 
\pages 47-157
\endref
\ref
\key 9
\by N.~Hitchin
\paper The geometry of three-forms in six and seven dimensions
\paperinfo DG/0010054
\yr 2000
\endref
\ref
\key 10
\by  D.D.Joyce
\paper Special Lagrangian 3-folds and integrable systems,
\paperinfo   
DG/0101249 
\endref    
\ref
\key 11
\by D.D.Joyce
\paper  Singularities of special Lagrangian fibrations and the SYZ Conjecture
\paperinfo  DG/0011179
\endref     
\ref 
\key 12
\by K.Kodaira
\book Complex manifolds and deformation of complex structures
\bookinfo Grundlehren der Mathematischen Wissenschaften,
\vol 283
\publ Springer-Verlag, New York-Berlin
\yr1986
\endref
\ref 
\key 13
\by H.B.~Lawson,~Jr and M.~Michelsohn 
\book Spin Geometry 
\publ Princeton University press 
\yr 1989
\endref
\ref 
\key 14
\by R.~McLean
\paper Deformations of calibrated submanifolds
\publ Comm. Anal. Geom. 
\vol 6 no.4 
\yr 1998
\pages 705--747
\endref
\ref
\key 15
\by
S.Salamon, 
\book Riemannian geometry and holonomy groups
\bookinfo Pitman Research Notes in Mathematics Series
\vol 201
\publ Longman, Harlow 
\yr 1989
\endref
\ref 
\key 16 
\by A.~ Strominger, S.T.~ Yau, E.~ Zaslow
\paper  Mirror Symmetry is T-Duality
\jour Nucl.Phys
\yr 1996
\vol B479
\pages 243-259
\endref
\ref
\key 17
\by G.Tian
\paper
 Smoothness of the universal deformation space of 
compact Calabi-Yau manifolds and its Petersson-Weil metric
\book Mathematical aspects of string theory
(ed. S.-T. Yau),
 \bookinfo Advanced
Series in Mathematical Physics
\vol 18
\publ World Scientific Publishing Co., Singapore
\yr 1987
\pages 629--646.
\endref
\ref \key 18
\by
A.N.Todorov, 
\paper
The Weil-Petersson geometry of the moduli space of 
 SU$(n\geq 3)$ (Calabi-Yau) manifolds. I 
\jour Comm. Math. Phys. 
\vol 126
\yr 1989
\pages 325--346
\endref
\endRefs
\enddocument